\documentclass[reqno]{amsart}
\usepackage{mathrsfs}
\usepackage{hyperref}

\AtBeginDocument{{\noindent\small
\emph{Electronic Journal of Differential Equations},
Vol. 2009(2009), No. 37, pp. 1--32.\newline
ISSN: 1072-6691. URL: http://ejde.math.txstate.edu or http://ejde.math.unt.edu
\newline ftp ejde.math.txstate.edu}
\thanks{\copyright 2009 Texas State University - San Marcos.}
\vspace{8mm}}

\begin{document}
\title[\hfilneg EJDE-2009/37\hfil Controllability, observability, realizability]
{Controllability, observability, realizability, and stability of
dynamic linear systems}

\author[J. M. Davis, I. A. Gravagne, B. J. Jackson, R. J. Marks II
\hfil EJDE-2009/37\hfilneg]
{John M. Davis, Ian A. Gravagne, Billy J. Jackson, Robert J. Marks II}

\address{John M. Davis \newline
Department of Mathematics\\
Baylor University\\
Waco, TX 76798, USA}
\email{John\_M\_Davis@baylor.edu}

\address{Ian A. Gravagne \newline
Department of Electrical and Computer Engineering\\
Baylor University\\
Waco, TX 76798, USA}
\email{Ian\_Gravagne@baylor.edu}

\address{Billy J. Jackson\newline
Department of Mathematics and Computer Science\\
Valdosta State University\\
Valdosta, GA 31698, USA}
\email{bjackson@valdosta.edu}

\address{Robert J. Marks II\newline
Department of Electrical and Computer Engineering\\
Baylor University\\
Waco, TX 76798, USA}
\email{Robert\_Marks@baylor.edu}

\thanks{Submitted January 23, 2009. Published March 3, 2009.}
\thanks{Supported by NSF Grants EHS\#0410685 and CMMI\#726996}
\keywords{Systems theory; time scale; controllability; observability;
 realizability; \hfill\break\indent
 Gramian; exponential stability; BIBO stability;
  generalized Laplace transform; convolution}

\begin{abstract}
 We develop a linear systems theory that coincides with the existing
 theories for continuous and discrete dynamical systems, but that
 also extends to linear systems defined on nonuniform time scales.
 The approach here is based on generalized Laplace transform methods
 (e.g. shifts and convolution) from the recent work \cite{DaGrJaMaRa}.
 We study controllability in terms of the controllability Gramian and
 various rank conditions (including Kalman's) in both the time invariant
 and time varying settings and compare the results. We explore
 observability in terms of both Gramian and rank conditions and
 establish related realizability results. We conclude by applying this
 systems theory to connect exponential and BIBO stability problems in
 this general setting. Numerous examples are included to show the
 utility of these results.
\end{abstract}

\maketitle
\numberwithin{equation}{section}
\newtheorem{theorem}{Theorem}[section]
\newtheorem{lemma}[theorem]{Lemma}
\newtheorem{definition}[theorem]{Definition}
\newtheorem{remark}[theorem]{Remark}
\newtheorem{example}[theorem]{Example}
\allowdisplaybreaks

\section{Introduction}

In this paper, our goal is to develop the foundation for a
comprehensive  linear systems theory which not only coincides with
the existing canonical systems theories in the continuous and
discrete cases, but also to extend those theories to dynamical
systems with nonuniform domains (e.g. the quantum time scale used
in quantum calculus \cite{kac}). We quickly see that the standard
arguments on $\mathbb{R}$ and $\mathbb{Z}$ do not go through when the graininess
of the underlying time scale is not uniform, but we overcome this
obstacle by taking an approach rooted in recent generalized
Laplace transform methods \cite{DaGrJaMaRa,Ja}. For those not
familiar with the rapidly expanding area of dynamic equations on
time scales, excellent references are \cite{bp1,bp2}.

We examine the foundational notions of controllability,
observability,  realizability, and stability commonly dealt with
in linear systems and control theory \cite{AnMi,cd,MiHoLi,rugh}.
Our focus here is how to generalize these concepts to the
nonuniform domain setting while at the same time preserving and
unifying the well-known bodies of knowledge on these subjects in
the continuous and discrete cases. This generalized framework has
already shown promising application to adaptive control regimes
\cite{ddgm,dgm}.\pagebreak

Throughout this work, we assume the following:
    \begin{itemize}
 \item $\mathbb{T}$ is a time scale that is unbounded above but with bounded
   graininess,
 \item $A(t) \in \mathbb{R}^{n \times n}, B(t) \in \mathbb{R}^{n \times m}, C(t)
   \in \mathbb{R}^{p \times n},$ and  $D(t) \in \mathbb{R}^{p \times m}$
   are all rd-continuous on $\mathbb{T}$,
 \item all systems in question are regressive.
    \end{itemize}
The third assumption implies that the matrix $I +\mu(t)A(t)$ is
invertible, and so on $\mathbb{Z}$, the transition matrix will always be
invertible. We are therefore justified in talking about
controllability rather  than reachability which is common
\cite{cd,eng,weiss} since the transition matrix in general need
not be invertible for $\mathbb{T}=\mathbb{Z}$.

In the following sections, we begin with the time varying case,
and then  proceed to treat the time invariant case. We will get
stronger (necessary and sufficient) results in the more
restrictive time invariant setting relative to the time varying
case (sufficient conditions). Although some of the statements
contained in this work can be found elsewhere
\cite{bkp,BaPiWy,fm}, in each of these cases proofs are either
missing, are restricted to time invariant systems, or are believed
to be in error \cite{fm} when $\mathbb{T}$ has nonuniform graininess.
Moreover---and very importantly---the methods used here are rooted
in Laplace transform techniques (shifts and convolution), and thus
are fundamentally different than the approaches taken elsewhere in
the literature. This tack overcomes the subtle problems that arise
in the arguments found in \cite{fm} when the graininess is
nonconstant.

\section{Controllability}

\subsection{Time Varying Case}
In linear systems theory, we say that a system is {\em controllable}
provided the solution of the relevant dynamical system
(discrete, continuous, or hybrid) can be driven to a specified final
state in finite time.  We make this precise now.

\begin{definition} \rm
Let $A(t)\in\mathbb{R}^{n \times n}$, $B(t)\in \mathbb{R}^{n \times m}$,
$C(t) \in \mathbb{R}^{p \times n}$, and $D(t)\in \mathbb{R}^{p \times m}$ all be
rd-continuous matrix functions defined on $\mathbb{T}$, with $p,m \leq n$.
The regressive linear system
    \begin{equation}\label{csystem}
    \begin{gathered}
    x^{\Delta}(t) =A(t)x(t)+B(t)u(t), \quad x(t_0)=x_0,\\
    y(t) =C(t)x(t) + D(t)u(t),
    \end{gathered}
    \end{equation}
is {\em controllable on $[t_0,t_f]$} if given any initial
state $x_0$ there exists a rd-continuous input $u(t)$ such
that the corresponding solution of the system satisfies
$x(t_f)=x_f$.
\end{definition}

Our first result establishes that a necessary and sufficient
condition for controllability of the linear system \eqref{csystem}
is the invertibility of an associated Gramian matrix.

\begin{theorem}[Controllability Gramian Condition]
\label{controlgramian}
The regressive linear system
    \begin{gather*}
    x^{\Delta}(t) = A(t)x(t)+B(t)u(t), \quad x(t_0)=x_0,\\
    y(t) = C(t)x(t) + D(t)u(t),
    \end{gather*}
is controllable on $[t_0,t_f]$ if and only if the $n\times n$
controllability Gramian matrix given by
    $$
    \mathscr{G}_{C}(t_0,t_f):=\int_{t_0}^{t_f}\Phi_{A}(t_0,\sigma(t))
B(t)B^{T}(t)\Phi^{T}_{A}(t_0,\sigma(t))\,\Delta t,
    $$
is invertible, where $\Phi_Z(t,t_0)$ is the transition matrix for the
system $X^\Delta(t) =Z(t)X(t)$, $X(t_0)=I$.
\end{theorem}

\begin{proof}
Suppose $\mathscr{G}_{C}(t_0,t_f)$ is invertible.  Then, given
$x_0$ and $x_f$, we can choose the input signal $u(t)$ as
    $$
    u(t)=-B^{T}(t)\Phi_{A}^{T}(t_{0},\sigma(t))\mathscr{G}_{C}^{-1}
(t_0,t_f)(x_0-\Phi_{A}(t_0,t_f)x_f),
\quad t \in [t_0,t_f),
    $$
and extend $u(t)$ continuously for all other
values of $t$.  The corresponding solution of the system at
$t=t_f$ can be written as
    \begin{align*}
    x(t_f)&=\Phi_{A}(t_f,t_0)x_0+\int_{t_0}^{t_f}\Phi_{A}(t_f,\sigma(t))B(t)u(t) \Delta t \\
    &=\Phi_{A}(t_f,t_0)x_0 \\
    &\quad -\int_{t_0}^{t_f}\Phi_{A}(t_f,\sigma(t))B(t)B^{T}(t)
 \Phi_{A}^{T}(t_0,\sigma(t)) \mathscr{G}_{C}^{-1}(t_0,t_f)(x_0
 -\Phi_{A}(t_0,t_f)x_f)\,\Delta t\\
    &=\Phi_{A}(t_f,t_0)x_0 -\Phi_{A}(t_f,t_0)\\
&\quad\times \int_{t_0}^{t_f}\Phi_{A}(t_0,\sigma(t))B(t)B^{T}(t)\Phi_{A}^{T}(t_0,\sigma(t))\,\Delta t \, \mathscr{G}_{C}^{-1}(t_0,t_f)(x_0-\Phi_{A}(t_0,t_f)x_f)\\
    &=\Phi_{A}(t_f,t_0)x_0-(\Phi_{A}(t_f,t_0)x_0-x_f)\\
    &=x_f,
    \end{align*}
so that the state equation is controllable on $[t_0,t_f]$.

For the converse, suppose that the state equation is controllable,
but for the sake of a contradiction, assume that the matrix
$\mathscr{G}_{C}(t_0,t_f)$ is not invertible.  If $\mathscr{G}_{C}(t_0,t_f)$ is not
invertible, then there exists a vector $x_a \neq 0$ such that
    \begin{align}\label{eq1}
    0=x_a^{T}\mathscr{G}_{C}(t_0,t_f)x_a&=\int_{t_0}^{t_f}x_a^{T}\Phi_{A}(t_0,\sigma(t))B(t)B^{T}(t)\Phi^{T}_{A}(t_0,\sigma(t))x_a\,\Delta t \notag\\
    &=\int_{t_0}^{t_f}\|x_a^{T}\Phi_{A}(t_0,\sigma(t))B(t)\|^2\,\Delta t,
    \end{align}
and hence
    \begin{equation}\label{eq2}
    x_a^T\Phi_{A}(t_0,\sigma(t))B(t)=0, \quad t \in [t_0,t_f).
    \end{equation}
However, the state equation is controllable on $[t_0,t_f]$, and so
choosing $x_0=x_a+\Phi_{A}(t_0,t_f)x_f$, there exists an input
signal $u_a(t)$ such that
    $$
    x_f=\Phi_{A}(t_f,t_0)x_0+\int_{t_0}^{t_f}\Phi_{A}(t_f,\sigma(t))B(t)u_{a}(t)\,\Delta t,
    $$
which is equivalent to the equation
    $$
    x_a=-\int_{t_0}^{t_f}\Phi_{A}(t_0,\sigma(t))B(t)u_{a}(t)\,\Delta t.
    $$
Multiplying through by $x_a^T$ and using \eqref{eq1} and \eqref{eq2}
yields $x_a^Tx_a=0$, a contradiction.  Thus, the matrix
$\mathscr{G}_{C}(t_0,t_f)$ is invertible.
\end{proof}

Since the controllability Gramian is symmetric and positive
semidefinite,  Theorem~\ref{controlgramian} can be interpreted as
saying \eqref{csystem} is controllable on $[t_0,t_f]$ if and only
if the Gramian is positive definite.  A system that is not
controllable on $[t_0,t_f]$ may become so when either $t_f$ is
increased or $t_0$ is decreased.  Likewise, a system that is
controllable on $[t_0,t_f]$ may become uncontrollable if $t_0$ is
increased and/or $t_f$ is decreased.

Although the preceding theorem is strong in theory, in practice it
is quite limited since computing the controllability Gramian
requires explicit knowledge of the transition matrix, but the
transition matrix for time varying problems is generally not known
and can be difficult to approximate in some cases. This
observation motivates the following definition and our next
theorem.

\begin{definition}\label{Kdef}
If $\mathbb{T}$ is a time scale such that $\mu$ is sufficiently
differentiable  with the indicated rd-continuous derivatives,
define the sequence of $n \times m$ matrix functions
    \begin{align*}
    K_{0}(t)&:=B(t),\\
    K_{j+1}(t)&:=(I+\mu(\sigma(t))A(\sigma(t)))^{-1}K_{j}^{\Delta}(t)
    -\Big[(I+\mu(\sigma(t))A(\sigma(t)))^{-1}(\mu^{\Delta}(t)A(\sigma(t))\\
    &\quad+\mu(t)A^{\Delta}(t))(I+\mu(t)A(t))^{-1}+A(t)(I+\mu(t)A(t))^{-1}
    \Big]K_{j}(t),
    \end{align*}
$j=0,1,2,\dots$.
\end{definition}

A straightforward induction proof shows that for all $t,s$, we have
    $$
    \frac{\partial^j}{\Delta s^j}\left[\Phi_{A}(\sigma(t),\sigma(s))B(s)
\right]=\Phi_{A}(\sigma(t),\sigma(s))K_{j}(s), \quad j=0,1,2,\dots
    $$
Evaluation at $s=t$ yields a relationship between these matrices
and those in Definition~\ref{Kdef}:
    $$
    K_{j}(t)=\frac{\partial^j}{\Delta s^j}
 \left[\Phi_{A}(\sigma(t),\sigma(s))B(s)\right]{\Big|}_{s=t}, \quad
j=0,1,2,\dots
    $$
This in turn leads to the following sufficient condition for controllability.

\begin{theorem}[Controllability Rank Theorem]
\label{kmatrix}
Suppose $q\in\mathbb{Z}^+$ such that, for $t \in
[t_0,t_f]$, $B(t)$ is $q$-times rd-continuously differentiable and
both of $\mu(t)$ and $A(t)$ are $(q-1)$-times rd-continuously
differentiable. Then the regressive linear system
    \begin{gather*}
    x^{\Delta}(t)  = A(t)x(t)+B(t) u(t), \quad x(t_0)=x_0, \\
    y(t)  =C(t) x(t) + D(t) u(t), \notag
    \end{gather*}
is controllable on $[t_0,t_f]$ if for some $t_c \in [t_0,t_f)$, we
have
    $$
    \mathop{\rm rank}
        \begin{bmatrix}
        K_0(t_c) & K_1(t_c) & \dots & K_q(t_c)
        \end{bmatrix}=n,
    $$
where
    $$
    K_j(t)=\frac{\partial^j}{\Delta s^j}\left[\Phi_{A}(\sigma(t),
\sigma(s))B(s)\right]\Big|_{s=t}, \quad j=0,1,\dots,q.
    $$
\end{theorem}

\begin{proof}
Suppose there is some $t_c \in [t_0,t_f)$ such that the
rank condition holds.  For the sake of a contradiction, suppose
that the state equation is not controllable on $[t_0,t_f]$.  Then
the controllability Gramian $\mathscr{G}_{C}(t_0,t_f)$ is not invertible
and, as
in the proof of Theorem~\ref{controlgramian}, there exists a
nonzero $n \times 1$ vector $x_a$ such that
    \begin{equation}\label{zerovec}
    x_a^T\Phi_{A}(t_0,\sigma(t))B(t)=0, \quad t \in [t_0,t_f).
    \end{equation}
If we choose the nonzero vector $x_b$ so that
$x_b=\Phi_{A}^T(t_0,\sigma(t_c))x_a$, then (\ref{zerovec}) yields
    $$
    x_b^T\Phi_{A}(\sigma(t_c),\sigma(t))B(t)=0, \quad t \in [t_0,t_f).
    $$
In particular, at $t=t_c$, we have $x_b^TK_0(t_c)=0$.
Differentiating \eqref{zerovec} with respect to $t$,
    $$
    x_b^T\Phi_{A}(\sigma(t_c),\sigma(t))K_1(t)=0, \quad t \in[t_0,t_f),
    $$
so that $x_b^TK_1(t_c)=0$.  In general,
    $$
    \frac{d^j}{\Delta t^j}\left[x_b^T\Phi_{A}^T(\sigma(t_c),\sigma(t))B(t)
 \right]\Big|_{t=t_c}=x_b^TK_j(t_c)=0, \quad j=0,1,\dots,q.
    $$
Thus,
    $$
    x_b^T
    \begin{bmatrix}
    K_0(t_c) & K_1(t_c) & \dots & K_q(t_c)
    \end{bmatrix}=0,
    $$
which contradicts the linear independence of the rows guaranteed
by the rank condition.  Hence, the equation is controllable on
$[t_0,t_f]$.
\end{proof}

When $\mathbb{T}=\mathbb{R}$, the collection of matrices $K_j(t)$ above is such
that each member is the $j$th derivative of the matrix
$\Phi_{A}(\sigma(t),\sigma(s))B(s)=\Phi_{A}(t,s)B(s)$. This
coincides with the literature in the continuous case (see, for
example, \cite{AnMi,cd,rugh}).  However, while still tractable, in
general the collection $K_j(t)$ is nontrivial to compute. The
mechanics are more involved even on $\mathbb{Z}$, which is still a very
``tame" time scale.  Therefore, the complications of extending the
usual  theory to the general time scales case are evident even at
this early juncture.

Furthermore, the preceding theorem shows that if the rank
condition  holds for some $q$ and some $t_c \in [t_0,t_f)$, then
the linear state equation is controllable on {\em any} interval
$[t_0,t_f]$ containing $t_c$.  This strong conclusion partly
explains why the condition is only a sufficient one.

\subsection{Time Invariant Case}
We now turn our attention to establishing results concerning the
controllability of regressive linear time invariant systems.  The
generalized Laplace transform presented in \cite{DaGrJaMaRa,Ja}
allows us to attack the problem in ways that simply are not
available in the time varying case.

First we recall a result from DaCunha in order to establish a
preliminary  technical lemma.

\begin{theorem}\cite{dac2}
\label{expsum} For the system $X^\Delta(t)=AX(t)$, $X(t_0)=I$,
there exist scalar functions $\{\gamma_0(t,t_0), \dots
,\gamma_{n-1}(t,t_0)\} \subset C_{{\rm rd}}^{\infty}(\mathbb{T},\mathbb{R})$ such
that the unique solution has representation
    $$
    e_{A}(t,t_0)=\sum_{i=0}^{n-1}A^i\gamma_{i}(t,t_0).
    $$
\end{theorem}

\begin{lemma}\label{linalg}
Let $A,B \in \mathbb{R}^{n \times n}$ and $u:=u_{x_0}(t_f,\sigma(s))
\in C_{\textup{rd}}(\mathbb{T},\mathbb{R}^{n \times 1})$. Then
    \begin{align}\label{inteq}
    \mathop{\rm span}{ \Big \{ } \int_{t_0}^{t_f}e_{A}(s,t_0)Bu_{x_0}(t_f,\sigma(s))
\Delta s{\Big \} }= \mathop{\rm span}\{B, AB, \dots, A^{n-1}B\}.
    \end{align}
\end{lemma}

\begin{proof}
Let $\{\gamma_k(t,t_0)\}_{k=0}^{n-1}$ be the collection of
functions that decompose the exponential matrix as guaranteed by
Theorem~\ref{expsum}. This collection forms a linearly independent
set since it can be taken as the solution set of an $n$-th order
system of linear dynamic equations. Apply the Gram-Schmidt process
to generate an orthonormal collection
$\{\hat{\gamma}_k(t,t_0)\}_{k=0}^{n-1}$.  The two collections are
related by
    \begin{align*}
    &\begin{bmatrix}
    \gamma_0(t,t_0) & \gamma_1(t,t_0) & \cdots &  \gamma_{n-1}(t,t_0)
    \end{bmatrix} \notag\\
    &=
    \begin{bmatrix}
    \hat{\gamma}_0(t,t_0) & \hat{\gamma}_1(t,t_0) & \cdots & \hat{\gamma}_{n-1}(t,t_0)
    \end{bmatrix}
    \begin{bmatrix}
    p_{11} & p_{12} & \cdots & p_{1n} \\
    0 & p_{22} & \cdots & p_{2n} \\
    \vdots & \vdots & \vdots & \vdots \\
    0 & 0 & \cdots & p_{nn}
    \end{bmatrix},
    \end{align*}
where the matrix on the right is the triangular matrix obtained from the
$QR$ factorization of the vector consisting of the functions
$\{\gamma_k(t,t_0)\}_{k=0}^{n-1}$ on the left.

Using the $QR$ factorization, we can write the matrix exponential as
    \begin{align*}
    e_{A}(t,t_0)&=\sum_{k=0}^{n-1}\gamma_{k}(t,t_0)A^k \\
    &=\sum_{k=0}^{n-1}
    \begin{bmatrix}
    \hat{\gamma_0}(t,t_0) & \hat{\gamma_1}(t_,t_0) & \cdots
    & \hat{\gamma}_{n-1}(t,t_0)
    \end{bmatrix}
    p_k A^k,
    \end{align*}
where $p_k$ is the $k$-th column vector of the triangular matrix $R$.
It is worth recalling here that the entries on the diagonal of this matrix
are norms of nonzero vectors and are thus positive.  That is, $p_{ii}>0$
for all $i$.

Rewriting the integral from \eqref{inteq},
    \begin{align*}
    &\int_{t_0}^{t_f}e_{A}(s,t_0)Bu_{x_0}(t_f,\sigma(s))\Delta s \\
    &=\int_{t_0}^{t_f}\sum_{k=0}^{n-1}\gamma_k(s,t_0)A^kBu_{x_0}(t_f,\sigma(s))\Delta s \\
    &=\sum_{k=0}^{n-1}A^kB\int_{t_0}^{t_f}\gamma_{k}(s,t_0)u_{x_0}(t_f,\sigma(s))\Delta s\\
    &=\sum_{k=0}^{n-1}A^kB\int_{t_0}^{t_f}
    \begin{bmatrix}
    \hat{\gamma}_{0}(s,t_0) & \hat{\gamma}_1(s,t_0) & \cdots & \hat{\gamma}_{n-1}(s,t_0)
    \end{bmatrix}
    p_ku_{x_0}(t_f,\sigma(s))\,\Delta s.
    \end{align*}
Let
    $$
    y_k=\int_{t_0}^{t_f}
    \begin{bmatrix}
    \hat{\gamma}_{0}(s,t_0) & \hat{\gamma}_1(s,t_0) & \cdots & \hat{\gamma}_{n-1}(s,t_0)
    \end{bmatrix}
    p_ku_{x_0}(t_f,\sigma(s))\,\Delta s,
    $$
$k=0, 1, \dots, n-1$.
We will show that $\mathop{\rm span}\{y_0,y_1,\dots, y_{n-1}\}=\mathbb{R}^n.$ That is,
there exists some $u \in C_{\text{rd}}(\mathbb{T},\mathbb{R}^{n \times 1})$ so
that for any arbitrary but fixed collection of vectors $\{z_0,z_1,
\dots, z_{n-1}\} \subset \mathbb{R}^{n \times 1}$, the system
    \begin{gather*}
    \int_{t_0}^{t_f}p_{11}\hat{\gamma}_0(s,t_0)u_{x_0}(t_f,\sigma(s))
\Delta s
:= z_0 =
        \begin{bmatrix}
        z_{00} \\
        z_{01} \\
        \vdots \\
        z_{0(n-1)}
        \end{bmatrix} \\
    \int_{t_0}^{t_f}\left(\hat{\gamma}_0(s,t_0)p_{12}
+\hat{\gamma}_1(s,t_0)p_{22}\right)u_{x_0}(t_f,\sigma(s))\Delta s
:= z_1 =
        \begin{bmatrix}
        z_{10} \\
        z_{11} \\
        \vdots \\
        z_{1(n-1)}
        \end{bmatrix} \\
     \vdots \\
    \int_{t_0}^{t_f}\left(\hat{\gamma}_0(s,t_0)p_{1n}+\dots
+\hat{\gamma}_{n-1}(s,t_0)p_{nn}\right)u_{x_0}(t_f,\sigma(s))\Delta
s := z_{n-1} =
        \begin{bmatrix}
        z_{(n-1)0} \\
        z_{(n-1)1} \\
        \vdots \\
        z_{(n-1)(n-1)}
        \end{bmatrix}
    \end{gather*}
has a solution.

To accomplish this, we use the fact that the collection
$\hat{\gamma_k}(s,t_0)$ is orthonormal and search for a solution of
the form
    $$
    u_{x_0}(t_f,\sigma(s))=(u_j)=\left(\sum_{i=0}^{n-1}\beta_{i}^j
\hat{\gamma_{i}}(s,t_0)\right).
    $$
Starting with $u_0$, the equations become
    \begin{gather*}
    \int_{t_0}^{t_f}\hat{\gamma_0}p_{11}\sum_{i=0}^{n-1}
 \beta_{i}^0\hat{\gamma_{i}}\,\Delta s=z_{00} \\
    \int_{t_0}^{t_f}\left(\hat{\gamma_0}p_{12}+\hat{\gamma_1}p_{22}\right)
\sum_{i=0}^{n-1}\beta_{i}^0\hat{\gamma_{i}}\,\Delta s=z_{10} \\
   \dots \\ \int_{t_0}^{t_f}\left(\hat{\gamma}_0 p_{1n}
+\hat{\gamma}_1 p_{2n}+\cdots+\hat{\gamma}_{n-1}p_{nn}\right)
\sum_{i=0}^{n-1}\beta_{i}^0\hat{\gamma}_{i}\,\Delta s=z_{(n-1)0}.
    \end{gather*}
Since the system $\hat{\gamma_{k}}$ is orthonormal, we can simplify the
equations above using the fact that the integral of cross terms
$\hat{\gamma_{i}}\hat{\gamma_{j}}$, $i \neq j$, is zero.
After doing so, the system becomes a lower triangular system that can
be solved by forward substitution.  (The observation that $p_{ii} \neq 0$
is crucial here, since this is exactly what allows us to solve the system.)
For example, the first equation becomes
    \[
    \int_{t_0}^{t_f}\hat{\gamma}_0 p_{11}\sum_{i=0}^{n-1}\beta_{i}^0\hat{\gamma}_{i}\,\Delta s=\int_{t_0}^{t_f}\hat{\gamma_0}^2\beta_{0}^0p_{11}\Delta s
    =\beta_{0}^0p_{11}
    =z_{01},
    \]
so that $\beta_{0}^0=\frac{z_{00}}{p_{11}}$.
Using this value for $\beta_0^0$ in the second equation,
    \begin{align*}
    \int_{t_0}^{t_f}\left(\hat{\gamma}_0 p_{12}+\hat{\gamma}_1 p_{22}\right) \sum_{i=0}^{n-1}\beta_{i}^0\hat{\gamma}_{i}\,\Delta s &= \int_{t_0}^{t_f}\left(\hat{\gamma}_{0}p_{12}+\hat{\gamma}_1 p_{22}\right) \left(\beta_{0}^0\hat{\gamma}_0 +\beta_{1}^0\hat{\gamma}_1 \right)\Delta
s \\
    & =\frac{p_{12}}{p_{11}}z_{01}+\beta_{1}^0 p_{22} \\
    & =z_{10},
    \end{align*}
so that
$\beta_1^0=\frac{1}{p_{22}}z_{11}-\frac{p_{12}}{p_{11}p_{22}}z_{01}$.
We can continue solving the system in like manner by using forward
substitutions to find $\beta_{j}^0$ for all $j=0,1, \dots, n-1$,
which will in turn yield
$u_0=\sum_{i=0}^{n-1}\beta_{i}^0\hat{\gamma}_i$.  Repeating this
process for $u_1,u_2, \dots, u_{n-1}$, we find the correct linear
combinations of $\hat{\gamma}_k$ to solve the system, and so the
claim follows.
\end{proof}

We are now in a position to establish the following analogue of the
Controllability Rank Theorem (Theorem~\ref{kmatrix}).

\begin{theorem}[Kalman Controllability Rank Condition]\label{kalmancont}
The time invariant regressive linear system
    \begin{gather*}
    x^{\Delta}(t)  =Ax(t)+B u(t), \quad x(t_0)=x_0, \\
    y(t)  =C x(t) + D u(t), \nonumber
    \end{gather*}
is controllable on $[t_0,t_f]$ if and only if the $n \times nm$
controllability matrix
    $$
    \left[B\ AB\ \cdots\ A^{n-1}B\right]
    $$
has rank $n$.
\end{theorem}

\begin{proof}
Suppose the system is controllable, but for the sake of a
contradiction that the rank condition fails.  Then there exists
an $n \times 1$ vector $x_a$ such that
    $$
    x_a^TA^kB=0, \quad k=0, \dots, n-1.
    $$
Now, there are two cases to consider: either
$x_a^Tx_f=0$ or $x_a^Tx_f \neq 0$.

Suppose $x_a^Tx_f \neq 0$. Then for any $t$, the
solution at time $t$ is given by
    \begin{align*}
    x(t) & =\int_{t_0}^{t}e_{A}(t,\sigma(s))Bu_{x_0}(s)\,\Delta s + e_{A}(t,t_0)x_0\\
    & =e_{A}(t,0) * Bu(t) +e_{A}(t,0)x_0 \\
    & =Bu(t)*e_{A}(t,0)+e_{A}(t,0)x_0 \\
    & =\int_{t_0}^{t}e_{A}(s,t_0)Bu_{x_0}(t,\sigma(s))\,\Delta s + e_{A}(t,t_0)x_0,
    \end{align*}
where we have written the solution as a (time scale) convolution and
appealed to the commutativity of the convolution \cite{DaGrJaMaRa,Ja}.
Choose initial state $x_0=By$, where $y$ is arbitrary. Then, again by
commutativity of the convolution and Theorem~\ref{expsum},
    \begin{align*}
    x_a^Tx(t) & =x_a^T\int_{t_0}^{t}e_{A}(s,t_0)Bu_{x_0}(t,\sigma(s))\,\Delta s +x_a^Te_{A}(t,t_0)x_0 \\
    &=\int_{t_0}^{t}\sum_{k=0}^{n-1}\gamma_{k}(s,t_0)x_a^TA^kBu_{x_0}(t,\sigma(s))\,\Delta s+\sum_{k=0}^{n-1}\gamma_k(t,t_0)x_a^TA^kBy \\
    & =0,
    \end{align*}
so that $x_a^Tx(t)=0$  for all $t$. This is a contradiction since
$x_a^Tx(t_f)=x_a^Tx_f \neq 0$.

Now suppose $x_a^Tx_f=0$.  This time, we choose initial state
$x_0=e_{A}^{-1}(t_f,t_0)x_a$. Similar to the equation above,
    \begin{align*}
    x_a^Tx(t) &= \int_{t_0}^{t}\sum_{k=0}^{n-1}
\gamma_{k}(s,t_0)x_a^TA^kBu_{x_0}(t,\sigma(s))\,
\Delta s+x_a^Te_{A}(t,t_0)e_{A}^{-1}(t_f,t_0)x_a \\
    &=x_a^Te_{A}(t,t_0)e_{A}^{-1}(t_f,t_0)x_a.
    \end{align*}
In particular, at $t=t_f$, $x_a^Tx(t_f)=\|x_a\|^2 \neq
0,$ another contradiction.

Thus in either case we arrive at a
contradiction, and so controllability implies the rank condition.

Conversely, suppose that the system is not controllable.
Then there exists an initial state $x_0 \in \mathbb{R}^{n \times 1}$ such that for
all input signals $u(t) \in \mathbb{R}^{m \times 1}$, we have $x(t_f)\neq x_f$.
Again, it follows from the commutativity of the convolution that
    \begin{align*}
    x_f \neq x(t_f) & =\int_{t_0}^{t_f}e_{A}(t_f,\sigma(s))Bu_{x_0}(s)
 \,\Delta s + e_{A}(t_f,t_0)x_0 \\
    & =\int_{t_0}^{t_f}e_{A}(s,t_0)Bu_{x_0}(t_f,\sigma(s))
 \,\Delta s+e_{A}(t_f,t_0) x_0 \\
    & =\int_{t_0}^{t_f}\sum_{k=0}^{n-1}\gamma_{k}(s,t_0)A^kBu_{x_0}
(t_f,\sigma(s))\,\Delta s +e_{A}(t_f,t_0)x_0.
    \end{align*}
In particular,
    $$
    \sum_{k=0}^{n-1}A^kB\int_{t_0}^{t_f}\gamma_{k}(s,t_0)u_{x_0}
 (t_f,\sigma(s))\,\Delta s \neq x_f-e_{A}(t_f,t_0)x_0.
    $$
Notice that the last equation holds if and only if there is no
linear combination of the matrices $A^kB$ for $k=0,1,\dots,n-1,$
which satisfies
    $$
    \sum_{k=0}^{n-1}A^kB\alpha_k = x_f-e_{A}(t_f,t_0)x_0.
    $$
The fact that there is no such linear combination follows from
Lemma~\ref{linalg} once we realize that an argument similar to the one
given in the proof of this result holds if $m<n$.  Thus, the matrix
    $$
    \left[B\ AB\ \cdots\ A^{n-1}B\right]
    $$
cannot have rank $n$, and so we have shown that if the matrix has rank $n$,
then it is controllable by contraposition.
\end{proof}

The preceding theorem is commonly called the Kalman Rank Condition
after R. E. Kalman who first proved it in 1960  for the cases
$\mathbb{T}=\mathbb{R}$ and $\mathbb{T}=\mathbb{Z}$ (see \cite{kalman2,kalman1}). Therefore our
analysis has unified the two cases, but we have also extended
these results to an arbitrary time scale with bounded graininess.
However, it is important to point out that the proof here is not
the one that Kalman gave, which is the one classically used for
$\mathbb{R}$ and $\mathbb{Z}$ (see \cite{rugh} for example).  In these two special
cases, an observation about the particular form of the matrix
exponential on $\mathbb{R}$ and $\mathbb{Z}$ (due to the uniform graininess)
allows one to arrive at the result in a more straightforward
manner. The general time scale case requires another argument
altogether as demonstrated above.

We now look at an example illustrating Theorem~\ref{kalmancont}.

\begin{example}\label{kalmancontex} \rm
Consider the system
    \begin{gather*}
    x^{\Delta}(t)  =
        \begin{bmatrix}
        -\frac{8}{45} & \frac{1}{30} \\
        -\frac{1}{45} & -\frac{1}{10}
        \end{bmatrix}
    x(t)+
        \begin{bmatrix}
        2 \\
        1
        \end{bmatrix}
    u(t), \quad x(0)=
        \begin{bmatrix}
        5 \\ 2
        \end{bmatrix}, \\
    y(t)  =
        \begin{bmatrix}
        3 & 4
        \end{bmatrix}
        x(t).
    \end{gather*}
It is straightforward to verify that
\[
    \mathop{\rm rank}
    \begin{bmatrix}
    B & AB
    \end{bmatrix}
    = \mathop{\rm rank}
    \begin{bmatrix}
    2 & -29/90 \\
    1 & -13/90
    \end{bmatrix}
    = 2,
\]
so that the state equation is controllable by Theorem~\ref{kalmancont}.
\end{example}

The next theorem establishes that there is a state variable change
in the time invariant case that demonstrates the ``controllable part"
of the state equation.

\begin{theorem}\label{Pmatrix}
Suppose the controllability matrix for the time invariant regressive
linear system
    \begin{gather*}
    x^{\Delta}(t)=Ax(t)+B u(t), \quad x(t_0)=x_0, \\
    y(t) =C x(t),
    \end{gather*}
satisfies
    $$
    \mathop{\rm rank}
    \begin{bmatrix}B & AB & \cdots & A^{n-1}B
    \end{bmatrix}=q,
    $$
where $0<q<n$.  Then there exists an invertible matrix $P$ such that
    $$
    P^{-1}AP=
    \begin{bmatrix}
    \hat{A}_{11} & \hat{A}_{12} \\
    0_{(n-q) \times q} & \hat{A}_{22}
    \end{bmatrix},
    \quad P^{-1}B=
    \begin{bmatrix}
    \hat{B}_{11}\\
    0_{(n-q) \times m}
    \end{bmatrix},
    $$
where $\hat{A}_{11}$ is $q \times q$, $\hat{B}_{11}$ is $q \times m$, and
    $$
    \mathop{\rm rank}
    \begin{bmatrix}
    \hat{B}_{11} & \hat{A}_{11}\hat{B}_{11} & \cdots & \hat{A}_{11}^{q-1}\hat{B}_{11}
    \end{bmatrix}
    =q.
    $$
\end{theorem}

\begin{proof}
We begin constructing $P$ by choosing $q$ linearly independent
columns $p_1$, $p_2$,\dots, $p_q$, from the controllability matrix
for the system.  Then choose $p_{q+1},\dots, p_{n}$ as $n \times
1$ vectors so that
    $$
    P=
    \begin{bmatrix}
    p_1 & \cdots & p_q & p_{q+1} & \cdots & p_{n}
    \end{bmatrix}
    $$
is invertible.  Define $G$ so that $PG=B$.  Writing the $j$-th
column of $B$ as a linear combination of the linearly independent
columns of $P$  given by $p_1,p_2,\dots, p_q$, we find that the
last $n-q$ entries of the $j$-th column of $G$ must be zero.  This
argument holds for $j=1,\dots,m$, and so $G=P^{-1}B$ does indeed
have the desired form.

Now set $F=P^{-1}AP$, yielding
    $$
    PF=
    \begin{bmatrix}
    Ap_1 & Ap_2 & \cdots & Ap_n
    \end{bmatrix}.
    $$
The column vectors $Ap_1,\dots,Ap_q$ can be written as linear
combinations of $p_1,\dots,p_n$ since each column of $A^kB,
k\geq0$ can be written as a linear combination of these vectors.
As for $G$ above, the first $q$ columns of $F$ must have zeros as
the last $n-q$ entries.  Thus, $P^{-1}AP$ has the desired form.
Multiply the rank-$q$ controllability matrix by $P^{-1}$ to obtain
    \begin{align*}
    P^{-1}
    \begin{bmatrix}
    B & AB & \cdots & A^{n-1}B
    \end{bmatrix}
    &=
    \begin{bmatrix}
    P^{-1}B & P^{-1}AB & \cdots & P^{-1}A^{n-1}B
    \end{bmatrix} \\
    &=
    \begin{bmatrix}
    G & FG & \cdots & F^{n-1}G
    \end{bmatrix} \\
    &=
    \begin{bmatrix}
    \hat{B}_{11} & \hat{A}_{11}\hat{B}_{11} & \cdots & \hat{A}_{11}^{n-1}\hat{B}_{11} \\
    0 & 0 & \dots & 0
    \end{bmatrix}.
    \end{align*}
Since the rank is preserved at each step, applying the Cayley-Hamilton
theorem gives
    $$
    \mathop{\rm rank}
    \begin{bmatrix}
    \hat{B}_{11} & \hat{A}_{11}\hat{B}_{11} & \cdots & \hat{A}_{11}^{q-1}\hat{B}_{11}
    \end{bmatrix}
    =q.
    $$
\end{proof}

Next, we use the preceding theorem to prove the following.

\begin{theorem}\label{Eigen}
The time invariant regressive linear system
    \begin{gather*}
    x^{\Delta}(t)=Ax(t)+B u(t), \quad x(t_0)=x_0, \\
    y(t)=C x(t),
    \end{gather*}
is controllable if and only if for every scalar $\lambda$ the only
complex $n \times 1$ vector $p$ satisfying $p^TA=\lambda p^T$, $p^TB=0$
is $p=0$.
\end{theorem}

\begin{proof}
For necessity, note that if there exists $p \neq 0$ and a complex
$\lambda$ such that the equation given is satisfied, then
    \begin{align*}
    p^T
    \begin{bmatrix}
    B & AB & \cdots & A^{n-1}B
    \end{bmatrix}
    & =
    \begin{bmatrix}
    p^TB & p^TAB & \cdots & p^TA^{n-1}B
    \end{bmatrix} \\
    & =
    \begin{bmatrix}
    p^TB & \lambda p^TB & \cdots & \lambda^{n-1}p^TB
    \end{bmatrix},
    \end{align*}
so that the $n$ rows rows of the controllability matrix are linearly
 dependent, and hence the system is not controllable.

For sufficiency, suppose that the state equation is not controllable.
Then by Theorem~\ref{Pmatrix}, there exists an invertible $P$ such that
\[
    P^{-1}AP=
    \begin{bmatrix}
    \hat{A}_{11} & \hat{A}_{12} \\ 0_{(n-q) \times q} & \hat{A}_{22}
    \end{bmatrix},
    \quad P^{-1}B=
    \begin{bmatrix}
    \hat{B}_{11}\\0_{(n-q) \times m}
    \end{bmatrix},
\]
with $0<q<n$.  Let $p^T=\begin{bmatrix}0_{1\times q}& p_q^T\end{bmatrix}P^{-1}$, where $p_q$ is a left eigenvector for $\hat{A}_{22}$.  Thus, for some complex scalar $\lambda$, $p_q^T\hat{A}_{22}=\lambda p_q^T$, $p_q\neq 0.$  Then $p \neq 0$, and
    \begin{gather*}
    p^TB  =
        \begin{bmatrix}
        0 & p_q^T
        \end{bmatrix}
        \begin{bmatrix}
        \hat{B}_{11} \\ 0
        \end{bmatrix}
    =0, \\
    p^TA  =
        \begin{bmatrix}
        0 & p_q^T
        \end{bmatrix}
        \begin{bmatrix}
        \hat{A}_{11} & \hat{A}_{12} \\0 & \hat{A}_{22}
        \end{bmatrix}
    P^{-1}=
        \begin{bmatrix}
        0 & \lambda p_q^T
        \end{bmatrix}
    P^{-1}=\lambda p^T.
    \end{gather*}
Thus, the claim follows.
\end{proof}

The interpretation of Theorem~\ref{Eigen} is that in a controllable
time invariant system, $A$ can have no left eigenvectors that are
orthogonal to the columns of $B$.  This fact can then be used to prove
the next theorem.

\begin{theorem}
The time invariant regressive linear system
    \begin{gather*}
    x^{\Delta}(t)=Ax(t)+B u(t), \quad x(t_0)=x_0, \\
    y(t) =C x(t),
    \end{gather*}
is controllable if and only if
    $\mathop{\rm rank} \begin{bmatrix}
    zI-A & B
    \end{bmatrix}  =n  $
for every complex scalar $z$.
\end{theorem}

\begin{proof}
By Theorem~\ref{Eigen}, the state equation is not controllable if and
only if there exists a nonzero complex $n \times 1$ vector $p$ and complex
scalar $\lambda$ such that
    $$
    p^T
    \begin{bmatrix}
    \lambda I -A & B
    \end{bmatrix}
    p \neq 0.
    $$
But this condition is equivalent to $ \mathop{\rm rank}\begin{bmatrix}\lambda
I -A & B\end{bmatrix}<n $.
\end{proof}

\section{Observability}

Next, we turn our attention to observability of linear systems.
As before, we treat the time varying case first followed by the time
invariant case.

\subsection{Time Varying Case}
In linear systems theory, when the term {\em observability} is
used, it refers to the effect that the state vector has on the
output of the state equation.  As such, the concept is unchanged
by considering simply the response of the system to zero input.
Motivated by this, we define the following.

\begin{definition} \rm
The regressive linear system
    \begin{gather*}
    x^{\Delta}(t)  =A(t)x(t), \quad x(t_0)=x_0, \\
    y(t)  =C(t)x(t),
    \end{gather*}
is {\em observable on $[t_0,t_f]$} if any initial state
$x(t_0)=x_0$ is uniquely determined by the corresponding
response $y(t)$ for $t \in [t_0,t_f)$.
\end{definition}

The notions of controllability and observability can be thought of
as duals of one another, and so any theorem that we obtain for
controllability should have an analogue in terms of observability.
Thus, we begin by formulating observability in terms of an associated
Gramian matrix.

\begin{theorem}[Observability Gramian Condition]
The regressive linear system
    \begin{gather*}
    x^{\Delta}(t)  =A(t)x(t), \quad x(t_0)=x_0, \\
    y(t) =C(t)x(t),
    \end{gather*}
is observable on $[t_0,t_f]$ if and
only if the $n \times n$ {\em observability Gramian matrix}
    $$
    \mathscr{G}_{O}(t_0,t_f):=\int_{t_0}^{t_f}\Phi_{A}^T(t,t_0)C^{T}(t)C(t)
\Phi_{A}(t,t_0)\,\Delta t,
    $$
is invertible.
\end{theorem}

\begin{proof}
If we multiply the solution expression
    $$
    y(t)=C(t)\Phi_{A}(t,t_0)x_0,
    $$
on both sides by $\Phi_{A}^T(t,t_0)C(t)$ and integrate, we obtain
    $$
    \int_{t_0}^{t_f}\Phi_{A}^T(t,t_0)C^T(t)y(t)\Delta t
=\mathscr{G}_{O}(t_0,t_f)x_0.
    $$
The left side of this equation is determined by $y(t)$ for
$t \in [t_0,t_f)$, and thus this equation is a linear algebraic equation
in $x_0$.  If $\mathscr{G}_{O}(t_0,t_f)$ is invertible, then $x_0$ is
uniquely determined.

Conversely, if $\mathscr{G}_{O}(t_0,t_f)$ is not invertible, then there
exists a nonzero vector $x_a$ such that $\mathscr{G}_{O}(t_0,t_f)x_a=0$.
But then $x_a^T\mathscr{G}_{O}(t_0,t_f)x_a=0$, so that
    $$
    C(t)\Phi_{A}(t,t_0)x_a=0, \quad t \in [t_0,t_f).
    $$
Thus, $x(t_0)=x_0+x_a$ yields the same zero-input response for the system
as $x(t_0)=x_0$, and so the system is not observable on $[t_0,t_f]$.
\end{proof}

The observability Gramian, like the controllability Gramian, is
symmetric  positive semidefinite.  It is positive definite if and
only if the state equation is observable.

Once again we see that the Gramian condition is not very practical
as it  requires explicit knowledge of the transition matrix.
Thus, we present a sufficient condition that is easier to check
for observability.  As before, observability and controllability
can be considered dual notions to one another, and as such, proofs
of corresponding results are often similar if not the same.  Any
missing observability proofs missing below simply indicates that
the controllability analogue should be consulted.

\begin{definition} \rm
If $\mathbb{T}$ is a time scale such that $\mu$ is sufficiently
differentiable  with the indicated rd-continuous derivatives,
define the sequence of $p \times n$ matrix functions
    \begin{align*}
    L_{0}(t) &:=C(t), \\
    L_{j+1}(t) &:=L_j(t)A(t)+L_j^{\Delta}(t)(I+\mu(t)A(t)), \quad
    j=0,1,2,\dots
    \end{align*}
\end{definition}

As in the case of controllability, an induction argument shows
that
$$
L_j(t)=\frac{\partial^j}{\Delta
t^j}\left[C(t)\Phi_{A}(t,s)\right]\Big|_{s=t}.
$$
With this, an argument similar to the one before shows the
following:

\begin{theorem}[Observability Rank Condition]
Suppose $q\in\mathbb{Z}^+$ is such that, for $t \in
[t_0,t_f]$, $C(t)$ is $q$-times rd-continuously differentiable and
both of $\mu(t)$ and $A(t)$ are $(q-1)$-times rd-continuously
differentiable. Then the regressive linear system
    \begin{gather*}
    x^{\Delta}(t)  =A(t)x(t), \quad x(t_0)=x_0, \\
    y(t)  =C(t) x(t),
    \end{gather*}
is observable on $[t_0,t_f]$ if for some $t_c \in [t_0,t_f)$, we
have
    $$
    \mathop{\rm rank}
    \begin{bmatrix}
    L_0(t_c) \\
    L_1(t_c) \\
    \vdots \\
    L_q(t_c)
    \end{bmatrix}
    =n,
    $$
where
    $$
    L_j(t)=\frac{\partial^j}{\Delta s^j}\left[C(t)\Phi_{A}(t,s)\right]
\Big|_{s=t}, \quad j=0,1,\dots,q.
    $$
\end{theorem}

\subsection{Time Invariant Case}
Like controllability, observability has equivalent conditions that become necessary and sufficient in the time invariant case.  We begin with a Kalman rank condition for observability.

\begin{theorem}[{Kalman Observability Rank Condition}]\label{kalmanob}
\ The time invariant regressive linear system
    \begin{gather*}
    x^{\Delta}(t)  =Ax(t), \quad x(t_0)=x_0, \\
    y(t)  =C x(t),
    \end{gather*}
is observable on $[t_0,t_f]$ if and only if the $nm \times n$
observability matrix
\[
    \begin{bmatrix}
    C  \\ CA \\ \vdots \\ CA^{n-1}
    \end{bmatrix}
\]
has rank $n$.
\end{theorem}

\begin{proof}
Again, we show that the rank condition fails if and only if the
observability Gramian is not invertible. Thus, suppose that the rank
condition fails.  Then there exists a nonzero $n \times 1$ vector
$x_a$ such that
    $$
    CA^kx_a=0, \quad k=0, \dots, n-1.
    $$
This implies, using Theorem~\ref{expsum}, that
\begin{align*}
    \mathscr{G}_{O}(t_0,t_f)x_a & =\int_{t_0}^{t_f}e_{A}^T(t,t_0)
C^{T}Ce_{A}(t,t_0)x_a\,\Delta t\\
    &=\int_{t_0}^{t_f}e_{A}^T(t,t_0)C^{T}\sum_{k=0}^{n-1}
\gamma_k(t,t_0)CA^kx_a\,\Delta t \\
    &=0,
\end{align*}
so that the Gramian is not invertible.

Conversely, suppose that the Gramian is not invertible.  Then
there exists nonzero $x_a$ such that
    $
    x_a^T\mathscr{G}_{O}(t_0,t_f)x_a=0.
    $
As argued before, this implies
    $$
    Ce_A(t,t_0)x_a=0, \quad t \in [t_0,t_f).
    $$
At $t=t_0$, we obtain $Cx_a=0$, and differentiating $k$ times and
evaluating the result at $t=t_0$ gives
    $CA^kx_a=0$, $k=0,\dots, n-1.$
Thus,
    $$
    \begin{bmatrix}
    C  \\ CA \\ \vdots \\ CA^{n-1}
    \end{bmatrix}
    x_a=0,
    $$
and the rank condition fails.
\end{proof}

The proof of the preceding result demonstrates an important point
about  controllability and observability in the arbitrary time
scale setting: namely, proofs of similar results for the two
notions are often similar, {\it but can sometimes be very
different}.  Indeed, comparing the proof of the Kalman condition
for controllability with the proof of the Kalman condition for
observability highlights this contrast.

The following example uses Theorem~\ref{kalmanob}.

\begin{example} \label{kalmanobex} \rm
Consider the system
    \begin{gather*}
    x^{\Delta}(t) =
        \begin{bmatrix}
        -\frac{8}{45} & \frac{1}{30} \\ -\frac{1}{45} & -\frac{1}{10}
        \end{bmatrix}
    x(t)+
        \begin{bmatrix}
        2 \\ 1
        \end{bmatrix}
    u(t), \quad x(0)=
        \begin{bmatrix}
        5 \\ 2
        \end{bmatrix}, \\
    y(t) =
        \begin{bmatrix}
        3 & 4
        \end{bmatrix}
    x(t).
    \end{gather*}
From Example \ref{kalmancontex}, recall that the system is
controllable.  We claim the system is also observable.  This
follows from Theorem~\ref{kalmanob} since
    $$
    \mathop{\rm rank}
        \begin{bmatrix}
        C \\ CA
        \end{bmatrix}
    =\mathop{\rm rank}
        \begin{bmatrix}
        3 & 4 \\ -\frac{28}{45} & -\frac{3}{10}
        \end{bmatrix}
    =2.
    $$
\end{example}

The following three theorems concerning observability have proofs that
mirror their controllability counterparts, and so will not be given here.

\begin{theorem} \label{obsvarchange}
Suppose the observability matrix for the time invariant regressive linear
system
    \begin{gather*}
    x^{\Delta}(t)=Ax(t)+B u(t), \quad x(t_0)=x_0, \\
    y(t)=C x(t),
    \end{gather*}
satisfies
    $$
    \mathop{\rm rank}
    \begin{bmatrix}
    C \\ CA \\ \vdots \\ CA^{n-1}
    \end{bmatrix}
    =\ell,
    $$
where $0<\ell<n$.  Then there exists an invertible $n \times n$ matrix
$Q$ such that
    $$
    Q^{-1}AQ=
    \begin{bmatrix}
    \hat{A}_{11} & 0 \\ \hat{A}_{21} & \hat{A}_{22}
    \end{bmatrix},
    \quad CQ=
    \begin{bmatrix}
    \hat{C}_{11} & 0
    \end{bmatrix},
    $$
where $\hat{A}_{11}$ is $\ell \times \ell$, $\hat{C}_{11}$ is
$p \times \ell$, and
    $$
    \mathop{\rm rank}
    \begin{bmatrix}
    \hat{C}_{11} \\\hat{C}_{11}\hat{A}_{11} \\
\vdots \\ \hat{C}_{11}\hat{A}_{11}^{\ell-1}
    \end{bmatrix}
    =\ell.
    $$
\end{theorem}

The state variable change Theorem~\ref{obsvarchange} is constructed
by choosing $n-\ell$ vectors in the nullspace of the observability matrix,
and preceding them by $\ell$ vectors that yield a set of $n$ linearly
independent vectors.

\begin{theorem} \label{obseval}
The time invariant regressive linear system
    \begin{gather*}
    x^{\Delta}(t)=Ax(t)+B u(t), \quad x(t_0)=x_0, \\
    y(t) =C x(t),
    \end{gather*}
is observable if and only if for every complex scalar $\lambda$, the
only complex $n \times 1$ vector $p$ that satisfies $Ap=\lambda p$,
$Cp=0$ is $p=0$.
\end{theorem}

Again, Theorem~\ref{obseval} can be restated as saying that in an
observable time invariant system, $A$ can have no right eigenvectors
that are orthogonal to the rows of $C$.

\begin{theorem}
The time invariant regressive linear system
    \begin{gather*}
    x^{\Delta}(t)=Ax(t)+B u(t), \quad x(t_0)=x_0, \\
    y(t) =C x(t),
    \end{gather*}
is observable if and only if
    $$
    \mathop{\rm rank}
    \begin{bmatrix}
    C \\z I -A
    \end{bmatrix}
    =n,
    $$
for every complex scalar $z$.
\end{theorem}

\section{Realizability}

In linear systems theory, the term {\em realizability} refers to the
ability to characterize a known output in terms of a linear system
with some input.  We now make this precise.

\begin{definition} \rm
The regressive linear system
    \begin{gather*}
    x^{\Delta}  =A(t)x(t)+B(t)u(t), \quad x(t_0)=0, \\
    y(t)  =C(t)x(t),
    \end{gather*}
of dimension $n$ is a {\em realization} of the weighting
pattern $G(t,\sigma(s))$ if
    $$
    G(t,\sigma(s))=C(t)\Phi_{A}(t,\sigma(s))B(s),
    $$
for all $t,s$. If a realization of this system exists, then the weighting
pattern is {\em realizable}. The system is a {\em
minimal realization} if no realization of $G(t,\sigma(s))$ with dimension
less than $n$ exists.
\end{definition}

Notice that for the system
    \begin{gather*}
    x^{\Delta}(t)  =A(t)x(t)+B(t)u(t), \quad x(t_0)=0, \\
    y(t)  =C(t)x(t) + D(t)u(t),
    \end{gather*}
the output signal $y(t)$ corresponding to a given input $u(t)$ and
weighting pattern $G(t,\sigma(s))=C(t)\Phi_{A}(t,\sigma(s))B(s)$
is given by
    $$
    y(t)=\int_{t_0}^{t}G(t,\sigma(s))u(s)\,\Delta s+D(t)u(t),
\quad t \geq t_0.
    $$

When there exists a realization of a particular weighting response
$G(t,\sigma(s)$, there will in fact exist many since a change of
state variables will leave the weighting pattern unchanged.  Also,
there can be many different realizations of the same weighting
pattern that all have different dimensions.  This is why we are
careful to distinguish between realizations and minimal
realizations in our definition.

We now give equivalent conditions for realizability: as before,
we begin with the time variant case and then proceed to the time
invariant case.

\subsection{Time Varying Case}
The next theorem gives a characterization of realizable systems in general.

\begin{theorem}[{Factorization of $G(t,\sigma(s))$}]
The weighting pattern $G(t,\sigma(s))$ is realizable if and only
if there exist a rd-continuous matrix $H(t)$ that is of dimension
$q \times n$ and a rd-continuous matrix $F(t)$ of dimension $n
\times r$ such that
    $$
    G(t,\sigma(s))=H(t)F(\sigma(s)), \quad \text{for all }t,s.
    $$
\end{theorem}

\begin{proof}
Suppose there exist matrices $H(t)$ and $F(t)$
with $G(t,\sigma(s))=H(t)F(\sigma(s))$.  Then the system
    \begin{gather*}
    x^{\Delta}(t)  =F(t)u(t), \\
    y(t)  =H(t)x(t),
    \end{gather*}
is a realization of $G(t,\sigma(s))$ since the transition matrix
of the zero system is the identity.

Conversely, suppose that $G(t,\sigma(s))$ is realizable.  We may assume
that the system
    \begin{gather*}
    x^{\Delta}(t)  =A(t)x(t)+B(t)u(t), \\
    y(t)  =C(t)x(t),
    \end{gather*}
is one such realization.  Since the system is regressive, we may write
    $$
    G(t,\sigma(s))=C(t)\Phi_{A}(t,\sigma(s))B(s)=C(t)\Phi_{A}(t,0)
\Phi_{A}(0,\sigma(s))B(s).
    $$
Choosing $H(t)=C(t)\Phi_{A}(t,0)$ and $F(t)=\Phi_{A}(0,\sigma(t))B(t)$,
the claim follows.
\end{proof}

Although the preceding theorem gives a basic condition for realization
of linear systems, often in practice it is not very useful because
writing the weighting pattern in its factored form can be very difficult.
Also, as the next example demonstrates, the realization given by the
factored form can often be undesirable for certain analyses.

\begin{example} \rm
Suppose $\mathbb{T}$ is a time scale with $0\leq \mu \leq 2$.  Under this assumption,
$-1/4 \in \mathcal{R}^+(\mathbb{T})$. Then the weighting pattern
    $$
    G(t,\sigma(s))=e_{-1/4}(t,\sigma(s)),
    $$
has the factorization
    $$
    G(t,\sigma(s))=e_{-1/4}(t,\sigma(s))=e_{-1/4}(t,0)e_{\ominus(-1/4)}(\sigma(s),0).
    $$
By the previous theorem, a one-dimensional realization of $G$ is
    \begin{gather*}
    x^{\Delta}(t)  =e_{\ominus(-1/4)}(t,0)u(t), \\
    y(t)  =e_{-1/4}(t,0)x(t).
    \end{gather*}
This state equation has an unbounded coefficient and is not uniformly exponentially
stable (note that $e_{\ominus (-1/4)}(t,0)=e_{1/(4-\mu)}(t,0)$ is unbounded
since $1/(4-\mu)>0$).
However, the one-dimensional realization of $G$ given by
    \begin{gather*}
    x^{\Delta}(t)  =-\frac{1}{4}x(t)+u(t), \\
    y(t)  =x(t),
    \end{gather*}
does have bounded coefficients and is uniformly exponentially stable.
\end{example}

Before examining minimal realizations, some remarks are in order.
 First, note that the inverse and $\sigma$ operators commute:
    \begin{align*}
    P^{-1}(\sigma(t)) & =P^{-1}(t)+\mu(t)(P^{-1}(t))^{\Delta} \\
    & =P^{-1}(t)+\mu(t)(-P(\sigma(t)))^{-1}P^{\Delta}(t)P^{-1}(t) \\
    & =P^{-1}(t)-(P(\sigma(t))^{-1}(P(\sigma(t))-P(t))P^{-1}(t) \\
    & =(P(\sigma(t)))^{-1}.
    \end{align*}

Second, it is possible to do a variable change on the system
    \begin{gather*}
    x^{\Delta}(t)  =A(t)x(t)+B(t)x(t), \\
    y(t)  =C(t) x(t),
    \end{gather*}
so that the coefficient of $x^{\Delta}(t)$ in the new system is zero,
while at the same time preserving realizability of the system under
the change of variables.

Indeed, set $z(t)=P^{-1}(t)x(t)$ and note that
$P(t)=\Phi_A(t,t_0)$ satisfies
    $$
    (P(\sigma(t)))^{-1}A(t)P(t)-(P(\sigma(t)))^{-1}P^{\Delta}(t)=0.
    $$
If we make this substitution, then the system becomes
    \begin{gather*}
    z^{\Delta}(t)  =P^{-1}(\sigma(t))B(t)u(t), \\
    y(t)  =C(t)P(t)z(t).
    \end{gather*}
Thus, in terms of realizability, we may assume without loss of
generality that $A(t) \equiv 0$ by changing the system to the form given
above.

It is important to know when a given realization is minimal. The
following theorem gives a necessary and sufficient condition for
this in terms of controllability and observability.

\begin{theorem}[{Characterization of Minimal Realizations}]
\label{minrealtv}
Suppose the regressive linear system
    \begin{gather*}
    x^{\Delta}(t)  =A(t)x(t)+B(t)x(t), \\
    y(t)  =C(t) x(t),
    \end{gather*}
is a realization of the weighting pattern $G(t,\sigma(s))$.  Then
this realization is minimal if and only if for some $t_0$ and
$t_f>t_0$ the state equation is both controllable and observable
on $[t_0,t_f]$.
\end{theorem}

\begin{proof}
As argued above, we may assume without loss of generality that
$A(t) \equiv 0$.  Suppose the $n$-dimensional realization given is
not minimal.  Then there is a lower dimension realization of
$G(t,\sigma(s))$ of the form
    \begin{align*}
    z^{\Delta}(t) & =R(t)u(t), \\
    y(t) & =S(t)z(t),
    \end{align*}
where $z(t)$ has dimension $n_z<n$.  Writing the weighting
pattern in terms of both realizations produces
$C(t)B(s)=S(t)R(s)$ for all $t,s$.  Thus,
    $$
    C^{T}(t)C(t)B(s)B^{T}(s)=C^{T}(t)S(t)R(s)B^{T}(s),
    $$
for all $t,s$. For any $t_0$ and $t_f>t_0$, it is possible to integrate this
expression with respect to $t$ and then with respect to $s$ to
obtain
    $$
    \mathscr{G}_{O}(t_0,t_f)\mathscr{G}_{C}(t_0,t_f)
    =\int_{t_0}^{t_f}C^{T}(t)S(t)\,\Delta t\int_{t_0}^{t_f}R(s)
B^{T}(s)\,\Delta s.
    $$
The right
hand side of this equation is the product of an $n \times n_z$
matrix and an $n_z \times n$ matrix, and as such, it cannot have
full rank since the dimension of the space spanned by the product is
at most $n_z<n$.  Therefore,
$\mathscr{G}_{O}(t_0,t_f)$ and $\mathscr{G}_{C}(t_0,t_f)$ cannot be
simultaneously invertible.
The argument is independent of the $t_0$ and $t_f$ chosen, and so
sufficiency is established.

Conversely, suppose that the given state equation is a minimal
realization of the weighting pattern $G(t,\sigma(s))$, with $A(t)
\equiv 0$.  We begin by showing that if either
    $$
    \mathscr{G}_{C}(t_0,t_f)=\int_{t_0}^{t_f}B(t)B^{T}(t)\,\Delta t,
    $$
or
    $$
    \mathscr{G}_{O}(t_0,t_f)=\int_{t_0}^{t_f}C^{T}(t)C(t)\,\Delta t,
    $$
is singular for all $t_0$ and $t_f$, then minimality is violated.
Thus, there exist intervals $[t_0^a,t_f^a]$ and $[t_0^b,t_f^b]$ such
that $\mathscr{G}_{C}(t_0^a,t_f^a)$ and $\mathscr{G}_{O}(t_0^b,t_f^b)$
are both invertible.  If we let $t_0={\rm min}\{t_0^a,t_0^b\}$ and
$t_f={\rm max}\{t_f^a,t_f^b\}$, then the positive definiteness of the
observability and controllability Gramians yields that both
$\mathscr{G}_{C}(t_0,t_f)$ and $\mathscr{G}_{O}(t_0,t_f)$ are invertible.

To show this, we begin by supposing that for every interval
$[t_0,t_f]$ the matrix $\mathscr{G}_{C}(t_0,t_f)$ is not invertible.
Then, given
$t_0$ and $t_f$ there exists an $n \times 1$ vector $x=x(t_0,t_f)$
such that
    $$
    0=x^T\mathscr{G}_{C}(t_0,t_f)x=\int_{t_0}^{t_f}B(t)B^{T}(t)x\,\Delta t.
    $$
Thus, $x^{T}B(t)=0$ for $t \in [t_0,t_f)$.

We claim that there exists at least one such vector $x$ that
is independent of $t_0$ and $t_f$.  To this end, note that if
$\mathbb{T}$ is unbounded from above and below, then for each positive
integer $k$ there exists an $n \times 1$ vector $x_k$ with
    $$
    \|x_k\|=1, \quad x_k^TB(t)=0, \quad t \in [-k,k].
    $$
Thus, $\{x_k\}_{k=1}^\infty$ is a bounded sequence of $n \times 1$
vectors and by the Bolzano-Wier\-strauss Theorem, it has a
convergent subsequence since $\mathbb{T}$ is closed.  We label this
convergent subsequence by $\{x_{k{_j}}\}_{j=1}^{\infty}$ and
denote its limit by $x_0=\lim_{j \to \infty}x_{k_{j}}.$  Note
$x_0^TB(t)=0$ for all $t$, since for any given time $t_a$, there
exists a positive integer $J_a$ such that $t_a \in [-k_j,k_j]$ for
all $j \geq J_a$, which in turn implies $x_{k_{j}}^TB(t_a)=0$ for
all $j \geq J_a$.  Hence, $x_0^T$ satisfies $x_0^TB(t_a)=0$.

Now let $P^{-1}$ be a constant, invertible, $n \times n$ matrix
with bottom row $x_0^T$. Using $P^{-1}$ as a change of state
variables gives another minimal realization of the weighting
pattern, with coefficient matrices
    $$
    P^{-1}B(t)=
    \begin{bmatrix}
    \hat{B}_1(t) \\ 0_{1 \times m}
    \end{bmatrix},
    \quad C(t)P=
    \begin{bmatrix}
    \hat{C}_1(t) & \hat{C}_2(t)
    \end{bmatrix},
    $$
where $\hat{B}_1(t)$ is $(n-1) \times m$, and $\hat{C}_1(t)$ is
$p \times (n-1)$.  Then a straightforward calculation shows
$G(t,\sigma(s))=\hat{C}_1(t)\hat{B}_1(\sigma(s))$ so that the linear
state equation
    \begin{gather*}
    z^{\Delta}(t)=\hat{B}_1(t)u(t), \\
    y(t)=\hat{C}_1(t)z(t),
    \end{gather*}
is a realization for $G(t,\sigma(s))$ of dimension $n-1$. This
contradicts the minimality of the original $n$-dimensional
realization.  Thus, there must be at least one $t_0^a$ and one
$t_f^a>t_0^a$ such that $\mathscr{G}_{C}(t_0^a,t_f^a)$ is
invertible.

A similar argument shows that there exists at least one $t_0^b$
and one $t_f^b>t_0^b$ such that  $\mathscr{G}_{O}(t_0^b,t_f^b)$ is
invertible. Taking $t_0=\min\{t_0^a,t_0^b\}$ and
$t_f=\max\{t_f^a,t_f^b\}$ shows that the minimal realization of
the state equation is both controllable and observable on
$[t_0,t_f]$.
\end{proof}

\subsection{Time Invariant Case}
We now restrict ourselves to the time invariant case and use a
Laplace  transform approach to establish our results. Instead of
considering the time-domain description of the input-output
behavior given by
    $$
    y(t)=\int_{0}^{t}G(t,\sigma(s))u(s) \,\Delta s,
    $$
we examine the corresponding behavior in the $z$-domain. Laplace
transforming the equation above and using the Convolution Theorem
\cite{DaGrJaMaRa} yields $Y(z)=G(z)U(z)$.  The question is:  given
a transfer function $G(z)$, when does there exist a time invariant
form of the state equation such that
    $$
    C(zI-A)^{-1}B=G(z),
    $$
and when is this realization minimal?

To answer this, we begin by characterizing time invariant
realizations. In what follows, a {\em strictly-proper rational
function} of $z$ is a rational function of $z$ such that the
degree of the numerator is strictly less than the degree of the
denominator.

\begin{theorem} \label{tireal}
The $ p \times q $ transfer function $G(z)$ admits a time
invariant realization of the regressive linear system
    \begin{gather*}
    x^{\Delta}(t)  =Ax(t)+Bu(t), \\
    y(t)  =C x(t),
    \end{gather*}
if and only if each entry of $G(z)$ is a strictly-proper rational
function of $z$.
\end{theorem}

\begin{proof}
If $G(z)$ has a time invariant realization, then $G$ has
the form $G(z)=C(zI-A)^{-1}B.$  We showed in \cite{DaGrJaMaRa,Ja} that for
each Laplace transformable function $f(t)$, $F(z) \to 0$ as $z \to
\infty$, which in turn implies that $F(z)$ is a strictly-proper
rational function in $z$.  Thus, the matrix $(zI-A)^{-1}$ is a
matrix of strictly-proper rational functions, and $G(z)$ is a
matrix of strictly-proper rational functions since linear
combinations of such functions are themselves strictly-proper and
rational.

Conversely, suppose that each entry $G_{ij}(z)$ in the matrix is
strictly-proper and rational.  Without loss of generality, we can
assume that each polynomial in the denominator is monic (i.e. has
leading coefficient of 1).  Suppose
    $$
    d(z)=z^{r}+d_{r-1}z^{r-1}+\cdots+d_0
    $$
is the least common
multiple of the polynomials in  denominators.  Then $d(z)G(z)$
can be decomposed as a polynomial in $z$ with $p \times q$
constant coefficient matrices, so that
    $$
    d(z)G(z)=P_{r-1}z^{r-1}+\cdots + P_1z+P_0.
    $$
We claim that the $qr$-dimensional matrices given by
    \begin{align*}
    A=
    \begin{bmatrix}
    0_q & I_q & \dots & 0_q \\
    0_q & 0_q & \dots & 0_q \\
    \vdots & \vdots & \vdots & \vdots \\
    0_q & 0_q & \dots & 0_q \\
    -d_0I_q & -d_1I_q & \dots & -d_{r-1}I_q
    \end{bmatrix},
    \quad
    B=
    \begin{bmatrix}
    0_q  \\ 0_q \\ \vdots \\ 0_q \\ I_q
    \end{bmatrix},
    \quad
    C=
    \begin{bmatrix}
    P_0 \\ P_1 \\ \vdots \\ P_{r-1}
    \end{bmatrix},
    \end{align*}
form a realization of $G(z)$.  To see this, let
    $$
    R(z)=(zI-A)^{-1}B,
    $$
and partition the $qr \times q$ matrix $R(z)$ into $r$ blocks
$R_1(z), R_2(z), \dots, R_{r}(z)$, each of size $q \times q$.
Multiplying $R(z)$ by $(zI-A)$ and writing the result in terms of
submatrices gives rise to the relations
    \begin{gather}    \label{r1}
    R_{i+1}(z)=zR_{i}, \quad i=1, \dots, r-1, \\
    \label{r2}
    zR_{r}(z)+d_0R_{1}(z)+d_1R_{2}(z)+\cdots+d_{r-1}R_{r}(z)=I_q.
    \end{gather}
Using (\ref{r1}) to rewrite (\ref{r2}) in terms of
$R_1(z)$ gives
    $$
    R_1(z)=\frac{1}{d(z)}I_q,
    $$
and thus from \eqref{r1} again, we have
    $$
    R(z)=\frac{1}{d(z)}
    \begin{bmatrix}
    I_q  \\ zI_q \\ \vdots \\ z^{r-1}I_q
    \end{bmatrix}.
    $$
Multiplying by $C$ yields
    $$
    C(zI-A)^{-1}B=\frac{1}{d(z)}\left(P_0+zP_1+\cdots+z^{r-1}P_{r-1}\right)
=G(z),
    $$
which is a realization of $G(z)$.
\end{proof}

The realizations that are minimal are characterized in the following theorem,
which is repeated here for completeness sake in this easier case of
Theorem~\ref{minrealtv}.

\begin{theorem} \label{tirealmin}
Suppose the time invariant regressive linear system
    \begin{gather*}
    x^{\Delta}(t)  =Ax(t)+Bx(t), \\
    y(t)  =C x(t),
    \end{gather*}
is a realization of the transfer function $G(z)$.  Then this state
equation is a minimal realization of $G(z)$ if and only if it is
both controllable and observable.
\end{theorem}

\begin{proof}
Suppose the state equation is a realization of $G(z)$ that is not
minimal.  Then there is a realization of $G(z)$ given by
    \begin{gather*}
    z^{\Delta}(t)  =P z(t) + Q z(t), \\
    y(t)  =R z(t),
    \end{gather*}
with dimension $n_z<n$.  Thus,
$$
Ce_{A}(t,0)B=Re_{P}(t,0)Q,\quad t \geq 0.
$$
Repeated differentiation with respect to $t$, followed by
evaluation at $t=0$ yields
    $$
    CA^kB=RF^kQ, \quad k=0,1,\dots
    $$
Rewriting this information in matrix form for $k=0,1,\dots,2n-2$,
we see
    $$
    \begin{bmatrix}
    CB & CAB & \cdots & CA^{n-1}B \\ \vdots & \vdots & \vdots & \vdots\\ CA^{n-1}B & CA^nB & \cdots & CA^{2n-2}B
    \end{bmatrix}
    =
    \begin{bmatrix}
    RQ & RPQ & \cdots & RP^{n-1}Q \\ \vdots & \vdots & \vdots & \vdots\\ RP^{n-1}Q & RP^nQ & \cdots & RP^{2n-2}Q
    \end{bmatrix},
    $$
which can be rewritten as
    \begin{align*}
    \begin{bmatrix}
    C \\ CA \\ \vdots \\ CA^{n-1}
    \end{bmatrix}
    \left[B \ AB \ \cdots \ A^{n-1}B\right]
    =\begin{bmatrix}
    R \\ RP \\ \vdots \\ RP^{n-1}
    \end{bmatrix}
    \left[Q \ PQ \ \cdots \ P^{n-1}Q\right].
    \end{align*}
However, since the right hand side of the equation is the product of an
$n_zp \times n_z$ and an $n_z \times n_zm$ matrix, the rank of the
product can be no greater than $n_z$.  Thus, $n_z<n$,
which implies that that the realization given in the statement
of the theorem cannot be both controllable and observable.
Therefore, by contraposition a controllable
and observable realization must be minimal.

Conversely, suppose the state equation given in the statement of
the theorem is a minimal realization that is not controllable.
Then there exists an $n \times 1$ vector $y \neq 0$ such that
    $$
    y^T \left[B \ AB \ \cdots \ A^{n-1}B\right]=0,
    $$
which implies $y^TA^{k}B=0$ for all $k \geq 0$
by the Cayley-Hamilton theorem.  For $P^{-1}$ an invertible $n
\times n$ matrix with bottom row $y^T$, then a variable change of
$z(t)=P^{-1}x(t)$ produces the state equations
    \begin{gather*}
    z^{\Delta}(t)  =\hat{A}z(t)+\hat{B}u(t), \\
    y(t)  =\hat{C}z(t),
    \end{gather*}
which is also an $n$-dimensional minimal realization of $G(z)$.
Partition the coefficient matrices of the state equation above as
    $$
    \hat{A}=P^{-1}AP=
    \begin{bmatrix}
    \hat{A}_{11} & \hat{A}_{12} \\ \hat{A}_{21} & \hat{A}_{22}
    \end{bmatrix},
    \quad
    \hat{B}=P^{-1}B=
    \begin{bmatrix}
    \hat{B}_1 \\ 0
    \end{bmatrix},
    \quad
    \hat{C}=CP=
    \begin{bmatrix}
    \hat{C}_1 & CA
    \end{bmatrix},
    $$
where $\hat{A}_{11}$ is $(n-1) \times (n-1)$, $\hat{B}_1$ is
$(n-1) \times 1$, and $\hat{C}_1$ is $1 \times (n-1)$.    From
these partitions, it follows from the construction of $P$ that
$\hat{A}\hat{B}=P^{-1}AB$ has the form
    $$
    \hat{A}\hat{B}=
    \begin{bmatrix}
    \hat{A}_{11}\hat{B}_1 \\ \hat{A}_{21}\hat{B}_1
    \end{bmatrix}
    =
    \begin{bmatrix}
    \hat{A}_{11}\hat{B}_1 \\ 0
    \end{bmatrix}.
    $$
Since the bottom row of $P^{-1}A^{k}B$ is zero for all $k \geq 0$,
    $$
    \hat{A}^k\hat{B}=
    \begin{bmatrix}
    \hat{A}_{11}^k\hat{B}_1 \\ 0
    \end{bmatrix},
    \quad k \geq 0.
    $$
But, $\hat{A}_{11}, \hat{B}_{1}, \hat{C}_1$ give an
$(n-1)$-dimensional realization of $G(z)$ since
    \begin{align*}
    \hat{C}e_{\hat{A}}(t,0)\hat{B} & =
    \begin{bmatrix}
    \hat{C}_1 & \hat{C}_2
    \end{bmatrix}
    \sum_{k=0}^{\infty}\hat{A}^k\hat{B}h_{k}(t,0) \\
    & =
    \begin{bmatrix}
    \hat{C}_1 & \hat{C}_2
    \end{bmatrix}
    \sum_{k=0}^{\infty}
    \begin{bmatrix}
    \hat{A}_{11}^k\hat{B}_1 \\
    0
    \end{bmatrix}
    h_{k}(t,0) \\
    & =\hat{C}_1e_{\hat{A}_{11}}(t,0)\hat{B}_1,
    \end{align*}
so that the state equation in the statement of the theorem
is in fact not minimal, a contradiction.  A
similar argument holds if the system is assumed not to be
observable.
 \end{proof}

We now illustrate Theorem~\ref{tireal} and Theorem~\ref{tirealmin}
with an example.

\begin{example}
Consider the transfer function
    $$
    G(z)=\frac{9(37+300z)}{5+75z+270z^2}.
    $$
$G(z)$ admits a time invariant realization by Theorem~\ref{tireal}
since $G(z)$ is a strictly-proper rational function of $z$.  The
form of $G(z)$ indicates that we should look for a 2-dimensional
realization with a single input and single output.  We can write
    $$
    G(z)= \begin{bmatrix}
        3 & 4
        \end{bmatrix}
    \left(zI-
        \begin{bmatrix}
        -\frac{8}{45} & \frac{1}{30} \\ -\frac{1}{45} & -\frac{1}{10}
        \end{bmatrix}
    \right)^{-1}
        \begin{bmatrix}
        2 \\ 1
        \end{bmatrix},
    $$
so that a time invariant realization of $G(z)$ is given by
    \begin{gather*}
    x^{\Delta}(t)  =
    \begin{bmatrix}
    -\frac{8}{45} & \frac{1}{30} \\ -\frac{1}{45} & -\frac{1}{10}
    \end{bmatrix}
    x(t) +
    \begin{bmatrix}
    2 \\ 1
    \end{bmatrix}
    u(t), \quad x(0)=x_0, \\
    y(t)  =
    \begin{bmatrix}
    3 & 4
    \end{bmatrix}
    x(t).
    \end{gather*}
We showed in Example~\ref{kalmancontex} that this realization is
in fact controllable, and we showed in Example~\ref{kalmanobex} that
it is also observable.  Thus, Theorem~\ref{tirealmin} guarantees
that this realization of $G(z)$ is minimal.
\end{example}

\section{Stability}

We complete our foray into linear systems theory by considering
stability. In  \cite{psw}, P\"otzsche, Siegmund, and Wirth deal
with exponential stability.  DaCunha also deals with this concept
under a different definition in \cite{dac1,dac3} and emphasizes
the time varying case. We begin by revisiting exponential
stability in the time invariant case and then proceed to another
notion of stability commonly used in linear systems theory.

\subsection{Exponential Stability in the Time Invariant Case}
We start this section by revisiting the notion of exponential stability.
We are interested in both the time invariant and
time varying cases separately since it is often possible to obtain stronger
results in the time invariant case.

We have already noted that if $A$ is constant, then
$\Phi_{A}(t,t_0)=e_{A}(t,t_0)$.  In what follows, we will consider
time invariant systems with $t_0=0$ in order to talk about the
Laplace transform.

DaCunha defines {\em uniform exponential stability} as follows.

\begin{definition}\cite{dac3} \rm
The regressive linear system
    $$
    x^{\Delta}=A(t)x(t), \quad x(t_0)=x_0,
    $$
is {\em uniformly exponentially stable} if there exist constants
$\gamma,\lambda>0$
with $-\lambda \in \mathcal{R}^+$ such that for any $t_0$ and
$x(t_0)$, the corresponding solution satisfies
    $$
    \|x(t)\| \leq \|x(t_0)\|\gamma e_{-\lambda}(t,t_0), \quad t \geq t_0.
    $$
\end{definition}

With this definition of exponential stability, we can prove the next theorem.

\begin{theorem} \label{expstabint}
The time invariant regressive linear system
    $$
    x^{\Delta}(t)=Ax(t), \quad x(t_0)=x_0,
    $$
is uniformly exponentially stable if and only if for some $\beta>0$,
    $$
    \int_{t_0}^{\infty}\|e_{A}(t,t_0)\| \,\Delta t \leq \beta.
    $$
\end{theorem}

\begin{proof}
For necessity, note that if the system is uniformly exponentially stable,
then by \cite[Theorem 3.2]{dac1},
    $$
    \int_{0}^{\infty}\|e_{A}(t,0)\|\,\Delta t \leq \int_{0}^{\infty}\gamma e_{-\lambda}(t,0) \,\Delta t=\frac{\gamma}{\lambda},
    $$
so that the claim follows.

For sufficiency, assume the integral condition holds but for the
sake of contradiction that the system is not exponentially stable.
Then, again by \cite[Theorem 3.2]{dac1}, for all $\lambda, \gamma
>0$ with $-\lambda \in \mathcal{R}^{+}$, we have
    $
    \|e_{A}(t,0)\|>\gamma e_{-\lambda}(t,0).
    $
Hence,
    $$
    \int_{0}^{\infty}\|e_{A}(t,0)\| \,\Delta t >\int_{0}^{\infty}\gamma e_{-\lambda}(t,0) \,\Delta t
    =\frac{\gamma}{-\lambda}e_{-\lambda}(t,0){\Big |}_{0}^{\infty}
    =\frac{\gamma}{\lambda}.
    $$
In particular, if we choose $\gamma > \beta \lambda$, then
    $$
    \int_{0}^{\infty}\|e_{A}(t,0)\| \,\Delta t > \frac{\beta \lambda}{\lambda}=\beta,
    $$
a contradiction.
\end{proof}

Now consider the system
    $$
    x^{\Delta}(t)=Ax(t), \quad x(0)=I.
    $$
Transforming this system yields
    $$
    X(z)=(zI-A)^{-1},
    $$
which is the transform of $e_{A}(t,0)$. This result is unique as
argued  in \cite{DaGrJaMaRa,Ja}.  Note that this matrix contains
only strictly-proper rational functions of $z$ since we have the
formula
    $$
    (zI-A)^{-1}=\frac{\mathop{\rm adj}(zI-A)}{\det(zI-A)}.
    $$
Specifically, $\det(zI-A)$ is an $n$th degree polynomial in $z$,
while each entry of $\mathop{\rm adj}(zI-A)$ is a polynomial of degree at most $n-1$.
 Suppose
    $$
    \det(zI-A)=(z-\lambda_1)^{\psi_1}(z-\lambda_2)^{\psi_2}\cdots(z-\lambda_m)^{\psi_m},
    $$
where $\lambda_1,\lambda_2,\dots,\lambda_n$ are the distinct
eigenvalues of the $n \times n$ matrix $A$, with corresponding
multiplicities $\psi_1,\psi_2,\dots,\psi_m$.  Decomposing
$(zI-A)^{-1}$ in terms of partial fractions gives
    $$
    (zI-A)^{-1}=\sum_{k=1}^{m} \sum_{j=1}^{\psi_k} W_{kj}
\frac{1}{(z-\lambda_k)^{j}},
    $$
where each $W_{kj}$ is an $n \times n$ matrix of partial fraction
expansion coefficients given by
    $$
    W_{kj}=\frac{1}{(\psi_k-j)!}\frac{d^{\psi_k-j}}{dz^{\psi_k-j}}
\left[(z-\lambda_k)^{\psi_k}(zI-A)^{-1}\right]\Big|_{z=\lambda_k}.
    $$
If we now take the inverse Laplace transform of $(zI-A)^{-1}$ in the
form given above, we obtain the representation
    \begin{equation}\label{rep}
    e_{A}(t,0)=\sum_{k=1}^{m} \sum_{j=1}^{\psi_k} W_{kj}\frac{f_{j-1}(\mu,\lambda_k)}{(j-1)!}e_{\lambda_{k}}(t,0),
    \end{equation}
where $f_{j}(\mu,\lambda_k)$ is the sequence of functions obtained
from the residue calculations of the $j$th derivative in the
inversion formula.  For example, the first few terms in the
sequence are
    \begin{align*}
    f_{0}(\mu,\lambda_k) & =1, \\
    f_{1}(\mu,\lambda_k) & =\int_{0}^{t}\frac{1}{1+\mu \lambda_k}\,\Delta \tau, \\
    f_{2}(\mu,\lambda_k) & =\Big(\int_{0}^{t}\frac{1}{1+\mu \lambda_k}\,\Delta
  \tau \Big)^2 - \int_{0}^{t}\frac{\mu}{(1+\mu \lambda_k)^2}\,\Delta \tau, \\
    f_{3}(\mu,\lambda_k) & =\Big(\int_{0}^{t}\frac{1}{1+\mu \lambda_k}\,\Delta
  \tau \Big)^3 - 3\int_{0}^{t}\frac{\mu}{(1+\mu \lambda_k)^2}\,\Delta \tau \int_{0}^{t}\frac{1}{1+\mu \lambda_k}\,\Delta \tau \\
    &\quad+\int_{0}^{t}\frac{2 \mu^2}{(1+\mu \lambda_k)^3}\,\Delta \tau, \\
    &\vdots
    \end{align*}
Notice that if $\mu$ is bounded, then each $f_{j}(\mu,\lambda_k)$
can be bounded by a ``regular" polynomial of degree $j$ in $t$, call it
$a_{j}(t)$.  That is, $f_{j}$ can be bounded by functions of the
form $a_j(t)=a_jt^j+a_{j-1}t^{j-1}\cdots+a_0$.  This observation
will play a key role in the next theorem.  P\"otzsche, Siegmund,
and Wirth do prove this result in \cite{psw}, but our proof differs
from theirs in that we use new transform results to obtain it, while
they use other techniques. Note, however, that in the theorem we do
use their definition of exponential stability rather than the one
given by DaCunha.  For completeness, we remind the reader by restating
their definition here.

\begin{definition}\cite{psw} \rm
\label{expstabpsw}
For $t,t_0 \in \mathbb{T}$ and $x_0 \in \mathbb{R}^n$, the system
    $$
    x^{\Delta}=A(t)x, \quad x(t_0)=x_0,
    $$
is {\em uniformly exponentially stable} if there exists a constant
$\alpha>0$ such that for every $t_0 \in \mathbb{T}$ there exists a
$K\geq 1$ with
    $$
    \|\Phi_{A}(t,t_0)\|\leq Ke^{-\alpha(t-t_0)}\text{ for }t \geq t_0,
    $$
with $K$ being chosen independently of $t_0$.
\end{definition}

Recall that DaCunha's definition of uniform exponential stability of
a system will imply that the system is uniformly exponential stable
if we use P\"otzsche, Siegmund, and Wirth's definition of the concept,
but the converse need not be true in general.  Thus, DaCunha's definition
is weaker in this sense.

\begin{theorem}[Spectral Characterization of Exponential Stability]
\label{expstab}
Let $\mathbb{T}$ be a time scale which is unbounded above but has bounded graininess. The time invariant regressive linear system
    $$
    x^{\Delta}(t)=Ax(t), \quad x(t_0)=x_0,
    $$
is uniformly exponentially stable
\textup{(}in the sense of Definition~\ref{expstabpsw}\textup{)}
if and only if $\mathop{\rm spec}(A) \subset\mathcal{S}(\mathbb{T})$,
the regressive set of exponential stability for $\mathbb{T}$, given by
    $$
    \mathcal{S}(\mathbb{T}):=\{\lambda\in\mathbb{C}:\limsup_{T\to\infty}
{1\over T-t_0}\int_{t_0}^T \lim_{s\searrow \mu(t)}
{\log|1+s\lambda|\over s}\,\Delta t<0\}.
    $$
\end{theorem}

\begin{proof}
Suppose the eigenvalue condition holds.  Then, appealing to
Theorem~\ref{expstabint} and writing the exponential in the
explicit form given above in terms of the distinct eigenvalues
$\lambda_1,\lambda_2,\dots,\lambda_m$, we obtain
    \begin{align*}
    \int_{0}^{\infty}\|e_{A}(t,0)\|\,\Delta t
& =\int_{0}^{\infty}\Big\|\sum_{k=1}^{m} \sum_{j=1}^{\psi_k} W_{kj}
\frac{f_{j-1}(\mu,\lambda_k)}{(j-1)!}e_{\lambda_{k}}(t,0)\Big\|\,\Delta t\\
    & \leq\sum_{k=1}^{m} \sum_{j=1}^{\psi_k} \|W_{kj}\|\int_{0}^{\infty}
\Big|\frac{f_{j-1}(\mu,\lambda_k)}{(j-1)!}e_{\lambda_k}(t,0)\Big| \,\Delta t \\
    & \leq\sum_{k=1}^{m} \sum_{j=1}^{\psi_k} \|W_{kj}\|\int_{0}^{\infty}\left|a_{j-1}(t)e_{\lambda_k}(t,0)\right| \,\Delta t \\
    & \leq\sum_{k=1}^{m} \sum_{j=1}^{\psi_k} \|W_{kj}\|\int_{0}^{\infty}a_{j-1}(t)e^{-\alpha t} \,\Delta t \\
    & \leq\sum_{k=1}^{m} \sum_{j=1}^{\psi_k} \|W_{kj}\|\int_{0}^{\infty}a_{j-1}(t)
    e^{-\alpha t} \,dt
     <\infty.
    \end{align*}
Note that the last three lines hold by appealing to Definition~\ref{expstabpsw}.  Thus, by Theorem~\ref{expstabint} the system is uniformly exponentially stable.

Now, for the sake of a contradiction, assume that the eigenvalue
condition fails.  Let $\lambda$ be an eigenvalue of $A$ with
associated eigenvector $v$, with $\lambda \notin \mathcal{S}(\mathbb{C})$.
Direct calculation shows that the solution of the system $x^{\Delta}=Ax$, $x(0)=v$, is given by $x(t)=e_{\lambda}(t,0)v$.  From \cite{psw}, if $\lambda \notin \mathcal{S}(\mathbb{C})$, then
    $
    \lim_{t\to\infty}e_{\lambda}(t,0) \neq 0,
    $
so that we arrive at a contradiction.
\end{proof}

\subsection{BIBO Stability in the Time Varying Case}
Besides exponential stability, the concept of bounded-input, bounded-output
stability is also a useful property for a system to have.
As its name suggests, the notion is one that compares the supremum of
the output signal with the supremum of the input signal.
Thus, we define the term as follows.

\begin{definition} \label{bibotv} \rm
The regressive linear system
    \begin{gather*}
    x^{\Delta}(t)  =A(t)x(t)+B(t)u(t), \quad x(t_0)=x_0, \\
    y(t) =C(t)x(t),
    \end{gather*}
is said to be {\em uniformly bounded-input, bounded-output} (BIBO)
{\em stable} if there exists a finite constant $\eta$ such that
for any $t_0$ and any input $u(t)$ the corresponding
zero-state response satisfies
    $$
    \sup_{t \geq t_0}\|y(t)\| \leq \eta \sup_{t \geq t_0}\|u(t)\|.
    $$
\end{definition}

Note that we use the word ``uniform" to stress that the same $\eta$
works for all $t_0$ and all input signals.

The following characterization of BIBO stability is useful.

\begin{theorem} \label{tvbibo}
The regressive linear system
    \begin{gather*}
    x^{\Delta}(t)  =A(t)x(t)+B(t)u(t), \quad x(t_0)=x_0, \\
    y(t) =C(t)x(t),
    \end{gather*}
is uniformly bounded-input, bounded-output stable if and only
if there exists a finite constant $\rho$ such that for all $t,\tau$ with
$t \geq \tau$,
    $$
    \int_{\tau}^{t}\|G(t,\sigma(s))\|\,\Delta s \leq \rho.
    $$
\end{theorem}

\begin{proof}
Assume such a $\rho$ exists.  Then for any $t_0$ and any input signal,
the corresponding zero-state response of the state equation satisfies
    $$
    \|y(t)\|=\left|\left|\int_{t_0}^{t}C(t)\Phi_{A}(t,\sigma(s))B(s)u(s) \, \Delta s\right|\right|\leq\int_{t_0}^{t}\|G(t,\sigma(s))\|\,\|u(s)\| \, \Delta s, \quad t\geq t_0.
    $$
Replacing $\|u(s)\|$ by its supremum over $s \geq t_0$, and using the
integral condition, we obtain
    $$
    \|y(t)\|\leq\sup_{t \geq t_0}\|u(t)\|\int_{t_0}^{t}\|G(t,\sigma(s))\|\,
\Delta s \leq\rho \sup_{t \geq t_0}\|u(t)\|, \quad t \geq t_0.
    $$
Thus, the system is BIBO stable.

Conversely, suppose the state equation is uniformly BIBO stable.
Then there exists a constant $\eta$ so that, in particular, the
zero-state response for any $t_0$ and any input signal such that
$\sup_{t \geq t_0}\|u(t)\| \leq 1$ satisfies $\sup_{t \geq
t_0}\|y(t)\| \leq \eta.$  For the sake of a contradiction, suppose
no finite $\rho$ exists that satisfies the integral condition.
Then for any given $\rho>0$, there exist $\tau_{\rho}$ and
$t_{\rho}>\tau_{\rho}$ such that
    $$
    \int_{\tau_{\rho}}^{t_{\rho}}\|G(t_{\rho},\sigma(s))\|\,\Delta s>\rho.
    $$
In particular, if $\rho=\eta$, this implies that there exist
$\tau_{\eta}$, with $t_{\eta}>\tau_{\eta}$, and indices $i,j$ such
that the $i,j$-entry of the impulse response satisfies
    $$
    \int_{\tau_{\eta}}^{t_{\eta}}|G_{ij}(t_{\eta},\sigma(s))|\,\Delta s>\eta.
    $$
With $t_0=\tau_{\eta}$ consider the $m \times 1$ input signal
$u(t)$  defined for $t \geq t_0$ as follows:  set $u(t)=0$ for
$t>t_{\eta}$, and for $t \in [t_0,t_{\eta}]$ set every component
of $u(t)$ to zero except for the $j$-th component given by the
piecewise continuous signal
    $$
    u_j(t):=
    \begin{cases}
    1, & G_{ij}(t_{\eta},\sigma(t))>0, \\
    0, & G_{ij}(t_{\eta},\sigma(t))=0,\\
    -1, & G_{ij}(t_{\eta},\sigma(t))<0.
    \end{cases}
    $$
This input signal satisfies $\|u(t)\| \leq 1$ for all $t \geq
t_0$, but because of the integral condition above, the $i$-th
component of  the corresponding zero-state response satisfies
    $$
    y_i(t_{\eta})=\int_{t_0}^{t_{\eta}}G_{ij}(t_{\eta},\sigma(s))u_j(s)\,\Delta s
    =\int_{t_0}^{t_{\eta}}|G_{ij}(t_{\eta},\sigma(s))|\,\Delta s \\
    >\eta.
    $$
Since $\|y(t_{\eta})\| \geq |y_i(t_{\eta})|$, we arrive at a
contradiction  that completes the proof.
\end{proof}

We now wish to give conditions under which the notions of
exponential  stability and BIBO stability are equivalent.  To this
end, we begin with the following.

\begin{theorem} \label{biboexpstab}
Suppose the regressive linear system
    \begin{gather*}
    x^{\Delta}(t)=A(t)x(t)+B(t)u(t), \quad x(t_0)=x_0, \\
    y(t)=C(t)x(t),
    \end{gather*}
is uniformly exponentially stable, and there exist constants
$\beta$ and $\gamma$ such that $\|B(t)\| \leq \beta$  and
$\|C(t)\| \leq \alpha$ for all $t$. Then the state equation is
also uniformly bounded-input, bounded-output stable.
\end{theorem}

\begin{proof}
Using the bound implied by uniform exponential stability, we have
    \begin{align*}
    \int_{\tau}^{t}\|G(t,\sigma(s))\|\,\Delta s
& \leq \int_{\tau}^{t}\|C(t)\| \, \|\Phi_{A}(t,\sigma(s))\| \,\|B(s)\|\,\Delta s \\
    & \leq\alpha \beta \int_{\tau}^{t} \|\Phi_{A}(t,\sigma(s))\|\,\Delta s \\
    & \leq\alpha \beta \int_{\tau}^{t}\gamma e_{-\lambda}(t,\sigma(s))\,\Delta s \\
    & \leq\frac{\alpha \beta \gamma}{\lambda} \int_{\tau}^{t} \frac{\lambda}{1-\mu(s) \lambda} e_{\lambda / (1- \mu \lambda)}(s,t)\, \Delta s\\
    & =\frac{\alpha \beta \gamma}{\lambda} \left(1-e_{-\lambda}(t, \tau)\right) \\
    & \leq\frac{\alpha \beta \gamma}{\lambda}.
    \end{align*}
By Theorem~\ref{tvbibo}, the state equation
is also bounded-input, bounded-output stable.
\end{proof}

The following example illustrates the use of Theorem~\ref{biboexpstab}.

\begin{example} \rm
Let $\mathbb{T}$ be a time scale with $0 \leq \mu <\frac{1}{2}$.
Consider the system
    \begin{gather*}
    x^{\Delta}(t)  =
    \begin{bmatrix}
    -2 & 1 \\ -1 &  -\sin(t)-2
    \end{bmatrix}
    x(t)+
    \begin{bmatrix}
    \cos(t) \\ \sin(t)
    \end{bmatrix}
    u(t), \quad x(t_0)=x_0,\\
    y(t)  =
    \begin{bmatrix}
    1 & e_{-1}(t,0)
    \end{bmatrix}x(t),
    \end{gather*}
where here, $\sin(t)$ and $\cos(t)$ are the usual trigonometric
functions and not their time scale counterparts.  DaCunha shows that
the system is uniformly exponentially stable by applying
\cite[Theorem~4.2]{dac1} (also found in \cite[Theorem 3.2]{dac3}) with
the choice $Q(t)=I$.  For $t \geq 0$, we have $\|B(t)\|=1$ and
    $
    \|C(t)\|=\sqrt{1+(e_{-1}(t,0))^2} \leq \sqrt{2},
    $
since $p=-1 \in \mathcal{R}^+$ from our assumption on $\mathbb{T}$.
By Theorem~\ref{biboexpstab}, the state equation is also uniformly
bounded-input, bounded-output stable.
\end{example}

For the converse of the previous theorem, it is known on $\mathbb{T}=\mathbb{R}$ and
$\mathbb{T}=\mathbb{Z}$ that stronger hypotheses than simply having the system be BIBO
stable are necessary to establish exponential stability
(see \cite{am,as,MiHoLi,rk}).  At present, we lack an analogue of
this result for an arbitrary time scale in the time varying system case.
We will see that the time invariant case does allow for the equivalence
of the two notions in the general time scale case under certain conditions.

\subsection{BIBO Stability in the Time Invariant Case}

In order to extend the definition of BIBO stability to the time
invariant case, we first need the following definitions.

\begin{definition}\cite{DaGrJaMaRa,Ja} \rm
Let $u(t)\in C_{\text{prd-e2}}(\mathbb{T},\mathbb{R})$. The {\it shift of $u(t)$ by
$\sigma(s)$}, denoted by $u(t,\sigma(s))$, is given by
    $$
    u(t,\xi)=\mathscr{L}^{-1}\{U(z)e_{\ominus z}^\sigma(s,0)\},
    $$
where $U(z):=\mathscr{L}\{u(t)\}(z)$, and $\mathscr{L}$, $\mathscr{L}^{-1}$ denote the generalized
Laplace transform and its inverse.
\end{definition}

\begin{definition}\cite{DaGrJaMaRa,Ja} \rm
For $f,g\in C_{\text{prd-e2}}(\mathbb{T},\mathbb{R})$, the {\it convolution} $f*g$
is given by
    $$
    (f*g)(t)=\int_0^t f(\tau)g(t,\sigma(\tau))\,\Delta\tau.
    $$
\end{definition}

\begin{definition} \label{biboti} \rm
For any shift $u(t,\sigma(s))$ of the transformable function $u(t)$, the time invariant system
    \begin{align*}
    x^{\Delta}(t) & =Ax(t)+Bu(t), \quad x(t_0)=x_0, \\
    y(t)& =Cx(t),
    \end{align*}
is said to be {\em uniformly bounded-input, bounded-output stable}
if there exists a finite constant $\eta$ such that
the corresponding zero-state response
satisfies
    $$
    \sup_{t \geq 0}\|y(t)\| \leq \eta \sup_{t \geq 0}\sup_{s \geq 0}\|u(t,\sigma(s))\|.
    $$
\end{definition}

Note that Definitions \ref{bibotv} and \ref{biboti} are different: one
deals with the time varying case and the other with the time invariant case.
 The modified definition in the time invariant case says that the output
stays bounded over all shifts of the input.

\begin{theorem}
The time invariant regressive linear system
    \begin{gather*}
    x^{\Delta}(t)=Ax(t)+Bu(t), \quad x(t_0)=x_0, \\
    y (t)=Cx(t),
    \end{gather*}
is bounded-input, bounded-output stable if and only if there exists
a finite $\beta>0$ such that
    $$
    \int_{0}^{\infty}\|G(t)\|\,\Delta t \leq \beta.
    $$
\end{theorem}

\begin{proof}
Suppose the claimed $\beta>0$ exists.  For any time $t$, we have
    $$
    y(t)=\int_{0}^{t}Ce_{A}(t,\sigma(s))Bu(s)\,\Delta s
=\int_{0}^{t}Ce_{A}(s,0)Bu(t,\sigma(s))\,\Delta s,
    $$
since $y(t)$ is a convolution. Hence,
    $$
    \begin{aligned}
    \|y(t)\| & \leq\|C\|\int_{0}^{t}\|e_{A}(s,0)\|\,\|B\|\,\sup_{0 \leq s \leq t}\|u(t,\sigma(s))\| \,\Delta s \\
    & \leq\|C\|\int_{0}^{\infty}\|e_{A}(s,0)\|\Delta s \,\|B\| \,\sup_{s \geq 0}\|u(t,\sigma(s))\|,
    \end{aligned}
    $$
which implies
    $$
    \sup_{t \geq 0}\|y(t)\| \leq \|C\|\int_{0}^{\infty}\|e_{A}(s,0)\|\,
\Delta s\,\|B\|\,\sup_{t \geq 0}\sup_{s \geq 0}\|u(t,\sigma(s))\|.
    $$
If we choose $\eta=\|C\| \, \beta \, \|B\|$, the claim follows.

Conversely, suppose that the system is bounded-input bounded-output stable,
but for the sake of a contradiction that the integral is unbounded.  Then,
    $$
    \sup_{t \geq 0}\|y(t)\| \leq \eta \sup_{t \geq 0}\sup_{s \geq 0}\|u(t,\sigma(s))\|,
    $$
and
    $$
    \int_{0}^{\infty}\|G(t)\|\Delta t >\beta \text{ for all }\beta>0.
    $$
In particular, there exist indices $i,j$ such that
    $$
    \int_{0}^{\infty}|G_{ij}(t)|\Delta t>\beta.
    $$
Choose $u(t,\sigma(s))$ in the following manner:
set $u_k(t,\sigma(s))=0$ for all $k \neq j$, and define $u_j(t,\sigma(s))$ by
    $$
    u_j(t,\sigma(s)):=
    \begin{cases}
    1,& \text{if }G_{ij}(s)>0, \\
    0, & \text{if }G_{ij}(s)=0, \\
    -1, & \text{if }G_{ij}(s)<0.
    \end{cases}
    $$
Choose $\beta>\eta>0$.  Note $\sup_{t \geq 0}\sup_{s \geq
0}\|u(t,\sigma(s)\|\leq 1,$ so $\sup_{t \geq 0}\|y(t)\| \leq
\eta$.  However,
    \begin{align*}
    \sup_{t \geq 0}\|y(t)\| & =\sup_{t \geq 0}
\Big\|\int_{0}^{t}G(s)u(t,\sigma(s))\,\Delta s\Big\| \\
    & =\sup_{t \geq 0}\Big\|\int_{0}^{t}G_j(s)\cdot u_{j}(s)\, \Delta s\Big\| \\
    & \geq\sup_{t \geq 0}\int_{0}^{t}|G_{ij}(s)|\, \Delta s \\
    & = \int_{0}^{\infty}|G_{ij}(s)|\, \Delta s\\
    &> \beta >\eta,
    \end{align*}
which is a contradiction.  Thus, the claim follows.
\end{proof}

The next theorem demonstrates the equivalence of exponential and BIBO
stability in the time invariant case.  Recall that this is a notion we
 currently lack in the time varying case.

\begin{theorem}[Equivalence of BIBO and Exponential Stability]
\label{tibiboexpeq}
Suppose the \\ time invariant regressive linear system
    \begin{gather*}
    x^{\Delta}(t)  =Ax(t)+B u(t), \quad x(t_0)=x_0, \\
    y(t)  =C x(t),
    \end{gather*}
is controllable and observable.  Then the system is uniformly
bounded-input, bounded output stable if and only if it is
exponentially stable.
\end{theorem}

\begin{proof}
If the system is exponentially stable, then
    $$
    \int_{0}^{\infty}\|Ce_{A}(t,0)B\|\,\Delta t \leq \|C\| \,
\|B\|\int_{0}^{\infty}\|e_{A}(t,0)\| \,\Delta t \leq \eta,
    $$
by Theorem~\ref{expstabint}.

Conversely, suppose the system is uniformly bounded-input, bounded
output stable.  Then
    $$
    \int_{0}^{\infty}\|Ce_{A}(t,0)B\| \,\Delta t<\infty,
    $$
which implies
    \begin{equation}
    \label{limitcondition}
    \lim_{t \to \infty}Ce_{A}(t,0)B=0.
    \end{equation}
Using the representation of the matrix exponential from
\eqref{rep}, we may write
\begin{equation}    \label{laplacerep}
    Ce_{A}(t,0)B =\sum_{k=1}^{m} \sum_{j=1}^{\psi_k}
    N_{kj}\frac{f_{j-1}(\mu,\lambda_k)}{(j-1)!}e_{\lambda_{k}}(t,0),
\end{equation}
 where the $\lambda_k$ are the distinct eigenvalues of $A$,
the $N_{kj}$ are constant matrices, and the $f_{j}(\mu,\lambda_k)$
are the terms from the residue calculations.  In this form,
    \begin{align*}\label{laplacederrep}
    &\frac{d}{\Delta t} Ce_{A}(t,0)B\\
    &=\sum_{k=1}^{m}\Big(N_{k1}\lambda_k +\sum_{j=2}^{\psi_k}
   \Big(\frac{f^{\Delta}_{j-1}(\mu,\lambda_k)(1+\mu(t)\lambda_k)}{(j-2)!}
   +\frac{\lambda_kf_{j-1}(\mu,\lambda_k)}{(j-1)!}\Big)\Big)
   e_{\lambda_k}(t,0).
    \end{align*}
If this function does not tend to
zero as $t\to\infty$, then using \eqref{laplacerep} and \eqref{limitcondition}, we arrive at a contradiction. Thus,
    $$
    \lim_{t \to \infty}\Big(\frac{d}{\Delta t}Ce_{A}(t,0)B\Big)
    =\lim_{t \to \infty}CAe_{A}(t,0)B=\lim_{t \to \infty}Ce_{A}(t,0)AB=0,
    $$
where the last equality holds by noting that if $A$ is constant,
then $A$ and $e_{A}(t,0)$ commute.  Similarly,
it can be shown that any order time derivative of the
exponential tends to zero as $t \to \infty$.  Thus,
    $$
    \lim_{t \to \infty}CA^ie_{A}(t,0)A^jB=0, \quad i,j=0,1,\dots
    $$
It follows that
    \begin{equation}\label{matrixlimit}
    \lim_{t \to \infty}
        \begin{bmatrix}
        C \\ CA \\ \vdots \\ CA^{n-1}
        \end{bmatrix}
    e_{A}(t,0)
        \begin{bmatrix}
        B & A B & \cdots & A^{n-1}B
        \end{bmatrix}
    =0.
    \end{equation}
But, the system is controllable and observable, and so we can form
invertible matrices $\mathscr{G}_{C}^a$ and $\mathscr{G}_{O}^a$ by
choosing $n$ independent columns of the controllability matrix and
$n$ independent rows of the observability matrix, respectively.
Then, by (\ref{matrixlimit}), $\lim_{t \to
\infty}\mathscr{G}_{O}^a e_{A}(t,0)\mathscr{G}_{C}^a=0$.  Hence,
$\lim_{t \to \infty}e_{A}(t,0)=0$ and exponential stability
follows from the arguments given in Theorem~\ref{expstab}.
\end{proof}

We  use  the preceding theorem in the following example.

\begin{example} \rm
Suppose $\mathbb{T}$ is a time scale with $0 \leq \mu \leq 4$.  The system
    \begin{gather*}
    x^{\Delta}(t)  =
    \begin{bmatrix}
    -\frac{8}{45} & \frac{1}{30} \\ -\frac{1}{45} & -\frac{1}{10}
    \end{bmatrix}
    x(t) +
    \begin{bmatrix}
    2 \\ 1
    \end{bmatrix}
    u(t), \quad x(0)=x_0, \\
    y(t)  =
    \begin{bmatrix}
    3 & 4
    \end{bmatrix}
    x(t).
    \end{gather*}
is controllable by Example~\ref{kalmancontex} and observable by
Example~\ref{kalmanobex}.  The eigenvalues of $A$ are
$\lambda_1=-\frac{1}{9}$ and $\lambda_2=-\frac{1}{6}$.
Note that the assumption on $\mathbb{T}$ implies
$\lambda_1, \lambda_2 \in \mathcal{S}(\mathbb{C})$, the stability region of $\mathbb{T}$.
Thus, by Theorem~\ref{expstab}, the system is exponentially stable.
Theorem~\ref{tibiboexpeq} then says that the system is also BIBO stable.
\end{example}

As we have seen, the Laplace transform can be a useful tool for analyzing
stability in the time invariant case.  With this in mind, we desire a
theorem that determines if a system is BIBO stable by examining its transfer
function.  The following theorem accomplishes this.

\begin{theorem}[{Transfer Function Characterization of BIBO Stability}]
The time invariant regressive linear system
    \begin{gather*}
    x^{\Delta}(t)=Ax(t)+B u(t), \quad x(t_0)=x_0, \\
    y(t)=C x(t),
    \end{gather*}
is bounded-input, bounded-output stable if and only if all poles of
the transfer function $G(z)=C(zI-A)^{-1}B$ are contained in
$\mathcal{S}(\mathbb{C})$.
\end{theorem}

\begin{proof}
If each entry of $G(z)$ has poles that lie in $\mathcal{S}(\mathbb{C})$,
then the partial fraction decomposition of $G(z)$ discussed earlier
shows that each entry of $G(t)$ is a sum of ``polynomial-times-exponential"
terms.  Since the exponentials will all have subscripts in the
stability region,
    \begin{equation}\label{gee}
    \int_{0}^{\infty}\|G(t)\|\,\Delta t<\infty,
    \end{equation}
and so the system is bounded-input, bounded-output stable.

Conversely, if \eqref{gee} holds,
then the exponential terms in any entry of $G(t)$ must have
subscripts  in the stability region by using a standard
contradiction argument.  Thus, every entry of $G(z)$ must have
poles that lie in the stability region.
\end{proof}

\end{document}